\documentclass[12pt]{article}
\usepackage{amssymb,amsmath,amsthm,mathrsfs}
\usepackage{float}

\usepackage{url}
\usepackage{bm}
\usepackage{graphicx}
\usepackage{algorithm}
\usepackage{subfigure}
\usepackage{algpseudocode}
\usepackage{subeqnarray}
\usepackage{cases}
\usepackage{authblk}
\usepackage{blindtext}
\usepackage{xcolor}
\usepackage{booktabs}
\usepackage{multirow}
\usepackage{appendix}

\newcommand{\la}{{\langle}}
\newcommand{\ra}{{\rangle}}

\newcommand{\I}{\mathcal{I}}

\newcommand{\bx}{\bm{x}}

\newcommand{\N}{\mathcal{N}}
\newcommand{\M}{\mathcal{M}}

\newcommand{\argmin}{\operatornamewithlimits{arg\,min}}
\newcommand{\argmax}{\operatornamewithlimits{arg\,max}}
\newtheorem{theorem}{Theorem}[section]
\newtheorem{corollary}[theorem]{Corollary}
\newtheorem{lemma}[theorem]{Lemma}
\newtheorem{assumption}[theorem]{Assumption}
\newtheorem{proposition}[theorem]{Proposition}
 \theoremstyle{definition}
\newtheorem{definition}[theorem]{Definition}
\theoremstyle{remark}
\newtheorem{remark}[theorem]{Remark}
\newtheorem{example}{Example}[section]
\providecommand{\keywords}[1]
{
	\small
	\textbf{\textit{Keywords---}} #1
}

\begin{document}
\date{}

\title{Minimization Over the Nonconvex Sparsity Constraint Using A Hybrid First-order method}

\author[1,2]{Xiangyu Yang}
\author[3]{Hao Wang}
\author[3]{Yichen Zhu}
\author[4]{Xiao Wang}
\affil[1]{School of Mathematics and Statistics, Henan University, Kaifeng 475000, China}
\affil[2]{Center for Applied Mathematics of Henan Province, Henan University, Zhengzhou, 450046, China}
\affil[3]{School of Information Science and Technology, ShanghaiTech University, Shanghai, 201210, China}
\affil[4]{Peng Cheng Laboratory, Shenzhen, 518066, China}
\affil[1,2]{\textit{yangxy@henu.edu.cn}}
\affil[3]{\textit{\{wanghao1,zhuych2022\}@shanghaitech.edu.cn}}
\affil[4]{\textit{wangx07@pcl.ac.cn}}
\maketitle

\begin{abstract}
We investigate a class of nonconvex optimization problems characterized by a feasible set consisting of level-bounded nonconvex regularizers, with a continuously differentiable objective. We propose a novel hybrid approach to tackle such structured problems within a first-order algorithmic framework by combining the Frank-Wolfe method and the gradient projection method. The Frank-Wolfe step is amenable to a closed-form solution, while the gradient projection step can be efficiently performed in a reduced subspace. A notable characteristic of our approach lies in its independence from introducing smoothing parameters, enabling efficient solutions to the original nonsmooth problems. We establish the global convergence of the proposed algorithm and show the  $O(1/\sqrt{k})$ convergence rate in terms of the optimality error for nonconvex objectives under reasonable assumptions. Numerical experiments underscore the practicality and efficiency of our proposed algorithm compared to existing cutting-edge methods. Furthermore, we highlight how the proposed algorithm contributes to the advancement of nonconvex regularizer-constrained optimization.
\end{abstract}


\keywords{Constrained optimization, Non-Lipschitz optimization, Frank-Wolfe variants, Weighted $\ell_1$-ball projection, Complexity analysis}

\section{Introduction}\label{Intro_sec}
Consider the following nonconvex sparsity-promoting constrained optimization problem
\begin{equation}\label{eq: main_Opt_General}
	\begin{aligned}
		\min_{\bm{x}\in\mathbb{R}^{n}} & \quad f(\bm{x})\\
		\textrm{s.t.}\                 & \quad \bm{x} \in \Omega := \{\bm{x}\mid\Phi(\bm{x}) := \sum_{i=1}^{n}\phi(\vert x_i\vert) \leq \gamma\},
	\end{aligned}
\end{equation}
where $\gamma \in (0,+\infty)$ is a user-specified parameter so that the feasible set $\Omega$ is nonempty and compact. The properties of $f$ and $\phi$ are assumed as follows.

\begin{assumption}\label{BasicAssum}
	(i) Function $f:\mathbb{R}^n \to \mathbb{R}$ is twice continuously differentiable but possibly nonconvex and has Lipschitz continuous gradient on $\Omega$ with modulus $L_f > 0$, namely, $\Vert\nabla f(\bm{x}) - \nabla f(\bm{y}) \Vert_2\leq L_f \Vert\bm{x} - \bm{y} \Vert_2, \forall \bm{x},\bm{y} \in \Omega$.
	\begin{itemize}
		\item[(ii)] Function $\phi:\mathbb{R}_{+} \to \mathbb{R}_{+}$ is continuous, strictly concave  with $\phi(0) =0$ and $\lim\limits_{t \to \infty} \phi(t) = \infty$. Moreover, it is differentiable on $\mathbb{R}_{+}\backslash\{0\}$ with $\phi'(t)\geq 0$ for any $t >0$, and is subdifferentiable at $0$. 
		
		\item[(iii)] The restriction of $\phi$ to $\mathbb{R}_{+}$ admits an inverse $\phi^{-1}: \mathbb{R}_{+} \to \mathbb{R}_{+}$. Hence, $\phi^{-1}$ is differentiable on $\mathbb{R}_{+}\backslash\{0\}$, and the right derivative $(\phi^{-1})'(0_{+}) $ also exists. 
	\end{itemize}
\end{assumption}

Problems of the form \eqref{eq: main_Opt_General} are motivated by sparsity-constrained optimization \cite{alcantara2022accelerated,bahmani2013greedy,bertsimas2016best,guo2021novel,lapucci2023unifying}. However, such kind of 
problems are NP-hard in general \cite{natarajan1995sparse}, posing challenges in both computation and analysis. To address this, researchers have explored nonconvex surrogates for the $\ell_0$ norm, such as Cappled $L_1$\cite{zhang2010analysis}, Smoothly
Clipped Absolute Deviation (SCAD)\cite{fan2001variable} and Minimax Concave Penalty (MCP)\cite{zhang2010nearly}, as alternative solutions, aiming to develop more effective and efficient algorithms. The past decade has witnessed great progress \cite{bian2017optimality,wang2021nonconvex}. The function $\phi$, determining the feasible region, can have various forms, including but not limited to those detailed in \ref{tab:Example_surrogate}.

\begin{table}[htbp]
	\caption{Concrete examples of constraint functions $\Phi$ and $(\kappa > 0 \textrm{ and }0<p<1)$.}
	\label{tab:Example_surrogate}
	\centering
	\begin{tabular}{cccc}
		\hline
		Regularizer & $\Phi$ &  $\phi'(\vert x_i\vert)$ &  $\phi^{-1}(\vert x_i\vert)$\\
		\hline
		\vspace{2mm}
		Exp\cite{bradley1998feature}  & $\sum_{i=1}^{n}1 - \exp^{-\kappa\vert x_i\vert}$ &   $\kappa \exp^{-\kappa\vert x_i\vert}$ &  $-\frac{1}{\kappa}\log(1-\vert x_i\vert)$ \\ 
		\hline
		\vspace{2mm}
		$\ell_p$ (quasi-)norm \cite{frank1993statistical} & $\sum_{i=1}^{n}\vert x_i\vert^{p}$ &  $p\vert x_i\vert^{p-1}$, $x_i \neq 0$ & $\vert x_i\vert^{1/p}$ \\ 
		\hline
		\vspace{2mm}
		Log \cite{lobo2007portfolio}& $\sum_{i=1}^{n}\log(1 + \kappa\vert x_i\vert)$ &   $\frac{\kappa}{1+\kappa\vert x_i\vert}$ &  $\frac{1}{\kappa}(\exp^{\vert x_i\vert}-1)$ \\ 
		\hline
		\vspace{2mm}
		Geman \cite{geman1995nonlinear}   & $\sum_{i=1}^{n}\frac{\vert x_i \vert}{\vert x_i \vert + \kappa}$ &    $\frac{\kappa}{(\vert x_i\vert +\kappa)^2}$ &  $\frac{\kappa \vert x_i\vert}{1-\vert x_i\vert}$ \\
		\hline 
		\vspace{2mm}
		Arctan\cite{candes2008enhancing}  & $\sum_{i=1}^{n}\arctan(\kappa\vert x_i\vert)$ &  $\frac{\kappa}{1 + (\kappa\vert x_i\vert)^2}$ & $\frac{1}{\kappa}\tan(\vert x_i\vert)$ \\ 
		\hline
	\end{tabular}%
\end{table} 

Motivated by the success of nonconvex $\ell_p$ norm, $p \in (0,1)$, in pursuing sparsity \cite{cheng2022interior,zhang2020smoothing}. We focus on the following $\ell_p$ ball-constrained optimization problem
\begin{equation}\label{eq: main_Opt}\tag{\text{$\mathscr{P}$}}
	\begin{aligned}
		\min_{\bm{x}\in\mathbb{R}^{n}} & \quad f(\bm{x})\\
		\textrm{s.t.}\                 & \quad \bm{x} \in \mathcal{B}_{\ell_p}  := \{\bm{x}\in \mathbb{R}^n \mid \Vert\bm{x}\Vert_{p}^{p}\leq \gamma\},
	\end{aligned}
\end{equation}
where $p \in (0,1)$. Problem \eqref{eq: main_Opt} arises in various fields of application.
\begin{example}[\textbf{Projection and compressive sensing}]
	For many machine learning and signal processing problems, if $f(\bm{x}) := \frac{1}{2}\Vert\bm{x}-\bm{y}\Vert_{2}^{2}$, where $\bm{y}\in\mathbb{R}^{n}$ is a given point, then \eqref{eq: main_Opt} corresponds to a Euclidean projection onto the nonconvex $\ell_{p}$ ball \cite{yang2022towards}. Furthermore, when a sensing matrix $\bm{A} \in \mathbb{R}^{m\times n}$ is available for use, \eqref{eq: main_Opt} could model the well-known $\ell_{p}$-constrained least-squares problems \cite{oymak2017sharp}, i.e., $f(\bm{x}) := \frac{1}{2}\Vert\bm{A}\bm{x} - \bm{y}\Vert_{2}^{2}$, in which $\bm{y}\in\mathbb{R}^{m}$ usually denotes observations polluted by noise.
\end{example}

\begin{example}[\textbf{Supervised sparse learning}] 
	With collected samples $\{(a_i,b_i)\}_{i=1}^{n} \subseteq \mathbb{R}^{n_{\bm{a}}}\times  \mathbb{R}^{n_{\bm{b}}}$ at hand, let loss function $\ell:\mathbb{R}^{n_{\bm{b}}} \times \mathbb{R}^{n_{\bm{b}}} \to \mathbb{R}$ (e.g., logistic loss and hinge loss for classification) be given.  One is typically interested in minimizing the empirical risk model, i.e., $f(\bm{x}) = \frac{1}{n}\sum_{i=1}^{n}\ell(h(a_i;\bm{x}),b_i)$ over $\mathcal{B}$ for training a useful learning model, in which $h(\cdot;\cdot)$ denotes a given prediction function.
\end{example}

\begin{example}[\textbf{Adversarial examples generation}]
	If $f(\bm{x}) := -g_{\text{adv}}(\bm{x} + \xi(\bm{x}))$, where $g_{\text{adv}}$ is a target classifier model that transforms the input $\bm{x}$ into a target class ``$\text{adv}$'' and $\xi(\bm{x})$ refers to an additive adversarial perturbation vector imposed on $\bm{x}$. Then \eqref{eq: main_Opt} becomes the pixel attacks problems in deep neural networks \cite{balda2020adversarial}. 
\end{example} 

Despite its widespread applications, \eqref{eq: main_Opt} presents challenges as it is nonsmooth, nonconvex, and non-Lipschitz when $p\in (0,1)$. As per our knowledge, limited progress has been made in solving \eqref{eq: main_Opt}, in contrast to the well-explored nonconvex $\ell_{p}$ (quasi)-norm penalized optimization problem \cite{candes2008enhancing,hu2017group,marjanovic2012l_q,xu2012l_}. Although the closely related unconstrained optimization model offers valuable insights for tackling $\ell_{p}$ ball-constrained problems \cite{yang2022towards}, algorithmic understanding for \eqref{eq: main_Opt} remains limited.   Therefore, the motivation of this paper is highlighted by developing a practical and efficient numerical algorithm for solving \eqref{eq: main_Opt} with desirable theoretical properties.

In this paper, we propose a novel hybrid first-order algorithm for solving \eqref{eq: main_Opt_General}, which is computationally efficient and boasts easy implementation. In particular, we focus on the representative $\ell_p$ ball-constrained optimization problem \eqref{eq: main_Opt}. The algorithm iterates between a Frank-Wolfe (FW) step and a gradient projection (GP) step. Specifically, at an iterate on the boundary of the $\ell_p$ ball, we replace the original nonconvex $\mathcal{B}_{\ell_p}$ by a well-constructed convex subset thereof and then define an appropriate quadratic function over such a projection-friendly subset. Consequently, for the next iterate, we solve a gradient projection subproblem in which the feasible set is a weighted $\ell_1$ ball formed by linearizing the $\ell_p$ ball at the current point. On the other hand, when the iterate is located in the interior of the feasible set, we obtain the search direction by solving an $\ell_p$ ball-constrained FW subproblem, which has a closed-form solution. If this search direction generates an infeasible point, we truncated it onto the boundary of the $\ell_p$ ball. This will lead to a gradient projection subproblem with a weighted $\ell_1$ ball in the next iteration. Our analysis shows that any cluster point of the iterate sequence is first-order optimal to \eqref{eq: main_Opt}. We further demonstrate a $O(1/\sqrt{k})$ convergence rate of the optimality error concerning the nonconvex objectives under reasonable assumptions. Moreover, the effectiveness and efficiency of the proposed algorithm are showcased by extensive numerical tests on $\ell_{p}$ ball projections and signal/image recovery tasks compared to state-of-the-art competitive algorithms.

\subsection{Related Literature and Our Contributions}
We mainly review recent progress in the literature for \eqref{eq: main_Opt}, involving the GP method and FW method that is closely related to our study. We then highlight the contributions of our proposed algorithm.

\textbf{Relation with GP method.} The gradient projection method described in \cite{jain2017non} for nonconvexly constrained problems can be applied for solving \eqref{eq: main_Opt}. However, its convergence is established only for strongly convex and strong smooth $f$ over the feasible region with a condition number smaller than $2$. Moreover, the projection onto a nonconvex set with the assumption that the \emph{global} optimal solution can be computed efficiently.  This is generally a relatively strong requirement since even pursuing a first-order stationary point of the $\ell_p$ ball projection subproblem is a nontrivial task \cite{yang2022towards}, not to mention global optimal solutions. In \cite{bahmani2013unifying}, the authors studied a GP method for minimizing the least-squares subject to the $\ell_{p}$ ball constraint. Based on the restricted isometry property conditions of the measurement matrix, the authors derived the linear convergence rate of the GP method. However, it is worth mentioning that their analysis also assumes that the \emph{global} optimal solution of the $\ell_p$-ball projection can be efficiently computed. As for the Euclidean projection onto the $\ell_{p}$ ball, i.e.,  $f(\bm{x}) = \tfrac{1}{2}\Vert\bm{x} - \bm{y}\Vert_{2}^{2}$,  the authors of \cite{yang2022towards} proposed an iteratively reweighted $\ell_1$ ball method, which was developed based on the majorization-minimization (MM) principle \cite{lange2016mm}. The key idea is first to add perturbation at each iteration to obtain a locally continuously differentiable $\ell_p$ norm function, and then to linearize the relaxed constraint to form a weighted $\ell_1$ ball. Since the projection is known on the boundary, the key of the algorithm in \cite{yang2022towards} is to drive the iterates to the boundary and avoid getting trapped around a local solution, which relies on the perturbation of the $\ell_p$ ball combined with a dynamic updating strategy to perturb. As for our case, the optimal solution is unclear to be on the boundary or in the interior of the feasible set. Therefore, we should turn to a different strategy that can allow the convergence to a stationary point in the interior of the $\ell_p$ ball. Moreover, our approach is free of any perturbations when performing projections onto the approximated feasible set. By noting that the relaxed nonconvex $\ell_{p}$ norm can be expressed as a difference-of-convex (DC) function \cite{fu1998penalized}, \cite{boob2020feasible} proposed a proximal point algorithm for solving \eqref{eq: main_Opt} with convergence and complexity analysis. To our knowledge, \cite{boob2020feasible} is the first work to consider the nonconvex sparsity-constrained model with \eqref{eq: main_Opt} as a particular case. In their work, strong convexity and Lipschitz smoothness are required for subproblems in order to seek a feasible subproblem solution and guarantee convergence; nevertheless, our algorithm is typically exempt from these conditions. Second, the authors presented a double-loop algorithm for each subproblem, which relies on efficient subproblem solvers. Last but not least, the updating rule for smoothing parameters should be carefully designed for convergence and some hyperparameters (e.g., initial radius) should be well-tuned in different problems. Meanwhile, the iteration complexity of their method ignores the impacts brought by smoothing parameters since they indeed considered a relaxed $\ell_{p}$ ball-constrained model.

\textbf{Relation with FW method.} The FW method \cite{frank1956algorithm} is an important class of first-order methods for solving convex set-constrained optimization problems. It is a projection-free method and is often considered as an alternative to the GP method. The FW method was originally applied to convex problems \cite{jaggi2013revisiting,lu2020generalized}. In the last decades, a number of research efforts have extended the FW method to the minimization of a nonconvex objective over convex sets \cite{kerdreux2021projection,wai2017decentralized}. As per our knowledge and as of very recently, \cite{zeng2024frankwolfetype} is the only work that considers \eqref{eq: main_Opt} via FW type methods for the optimization involving nonconvex feasible set in our case. In order to find a feasible direction, they considered implementing a linear-optimization oracle, in which the concerned constrained set is formed by employing affine minorization to approximate the second DC component at each iteration. They established the global convergence analysis that iterations cluster in first-order stationary points. Despite not being stated explicitly, its numerical performance may be impacted by the updating rule for smoothing parameter in $\ell_{p}$ norm at each iteration when taking into account the DC expression of the $\ell_{p}$ norm.

Given these developments, the following summarizes our contributions:
\begin{itemize}
	\item[1.] \textbf{A novel approach under hybrid first-order algorithmic framework with efficient subproblem solutions.} We take advantage of the FW method and GP method for solving \eqref{eq: main_Opt}. FW method enables the search for interior iterates, while the GP method intends for boundary iterates. Both are available for efficient subproblem solutions. It is also important to note that, for the GP subproblem, our algorithm is free of introducing a smoothing parameter, thereby only performing projections onto a weighted $\ell_1$ ball in a reduced subspace. This underscores the difference between our method and those of \cite{boob2020feasible,yang2022towards,zeng2024frankwolfetype}.
	
	\item[2.] \textbf{Global convergence and worst-case complexity analysis.} We prove that the iterates globally converge to stationary points of \eqref{eq: main_Opt}. Additionally, we establish that the worst-case iteration complexity to attain an $\epsilon$-optimal solution is in order $O(1/\epsilon^2)$. While our complexity result is in the same order as that of \cite{boob2020feasible} for the deterministic nonconvex and nonsmooth case, our subproblem resolution is much more efficient, as evidenced by our numerical comparison.
	
	\item[3.]  \textbf{The effectiveness and efficiency of the proposed algorithm.} We demonstrate the effectiveness and efficiency of the proposed algorithm through a set of numerical experiments using both synthetic and real-world data. Our algorithm is generally superior to other comparative algorithms, including \cite{boob2020feasible,yang2022towards,zeng2024frankwolfetype} on $\ell_{p}$ ball projection problems. In particular, our algorithm could handle the cases when $p$ is stringent (e.g., $p < 0.5$), in which other algorithms are not comparable. In addition, our proposed algorithm can solve the signal and image recovery tasks effectively and efficiently.
\end{itemize}

\subsection{Organization}
The paper is structured as follows.  We describe the proposed algorithm in \S\ref{Algo_sec}. The global convergence and complexity are analyzed in \S\ref{Conver_sec}.  We present the results of numerical experiments in \S\ref{Exp_sec}. Conclusions are provided in \S\ref{Con_sec}. 

\subsection{Notation and Preliminaries}
In this subsection, we set the notation and provide some necessary preliminaries used throughout this paper. Specifically, let $\mathbb{R}^n$ denote the real $n$-dimensional Euclidean space with the standard inner product $\la\cdot,\cdot\ra : \mathbb{R}^{n} \times \mathbb{R}^{n} \to \mathbb{R}$.  Let $\mathbb{R}^n_+$ represent the non-negative orthant in $\mathbb{R}^n$, and $\mathbb{R}^{n}_{++}$ denote the interior of $\mathbb{R}^{n}_{+}$. Correspondingly, $\mathbb{R}^n_-$ refers to the non-positive orthant in $\mathbb{R}^n$. We use $[n]$ to refer the index set $\{1,\ldots,n\}$. For $\bm{x}\in\mathbb{R}^{n}$, we denote by $x_i$
the $i$th entry of $\bm{x}$, and define $\vert \bm{x} \vert = (\vert x_1 \vert, \ldots, \vert x_n \vert)^{T}$. For any $\bm{x},\bm{y}\in\mathbb{R}^{n}$, $\bm{x} \geq \bm{y}$ means $x_{i} \geq y_{i}$, $i \in [n]$, and we use $\circ$ to denote the entrywise product between $\bm{x}, \bm{y}$, i.e., $(\bm{x}\circ\bm{y})_i = x_iy_i$.
The signum function of  $\bm{x}\in\mathbb{R}^{n}$ is represented by $\text{sgn}(\bm{x}) = (\text{sgn}(x_1),\ldots, \text{sgn}(x_n))^{T}$, where $\text{sgn}(x) = 1$ if $x > 0$, $\text{sgn}(x) = 0$ if $x = 0$ and $\text{sgn}(x) = -1$ if $x < 0$. 

Let $\I(\bm{x}) := \{i \in [n] \mid x_i \neq 0\}$ and $\mathcal{A}(\bm{x}) := \{i \in [n] \mid x_i = 0\}$ denote the  inactive set and active set of $\bm{x}$, respectively. For a set $\mathcal{C} \subseteq \mathbb{R}^{n}$, the indicator function of $\mathcal{C}$ is given by $\delta_{\mathcal{C}}(\bm{x}) = 0, \forall \bm{x} \in \mathcal{C}$ and $\delta_{\mathcal{C}}(\bm{x}) = +\infty, \forall \bm{x} \notin \mathcal{C}$.
The interior and boundary of $\mathcal{C} \subseteq \mathbb{R}^{n}$ are denoted as $\textbf{int}\ \mathcal{C}$ and $\textbf{bd}\ \mathcal{C}$, respectively. We use \textrm{conv}($\mathcal{C}$) to denote the \textit{convex hull} of a set $\mathcal{C} \subset \mathbb{R}^n$, which is explicitly defined as $\textrm{conv}(\mathcal{C}) := \{\sum_{i=1}^{n}\theta_ix_i \mid x_i \in \mathcal{C}, \theta_i \geq 0, \sum_{i=1}^{n}\theta_i =1, i \in [n]\}$. We follow from \cite[Eq. 3(11)]{rockafellar2009variational} to define the \textit{positive hull} of a set $\mathcal{C} \subseteq \mathbb{R}^{n}$ as $\textrm{pos}\ \mathcal{C} = \{0\} \cup \{t\bm{x}\mid \bm{x} \in \mathcal{C}, \ t >0\}$.
Additionaly,  the Gaussion distribution with mean $\mu$ and standard deviation $\sigma$ is represented by $\mathcal{N}(\mu,\sigma^{2})$.

We next recall the definitions of subdifferentials and normal cones that are standard tools in variational analysis for developing optimality conditions of \eqref{eq: main_Opt}, adapting from \cite[Definition 8.3]{rockafellar2009variational} and \cite[Definition 6.3]{rockafellar2009variational}.
\begin{definition}\label{Def_Subdifferential}
	Consider a proper function $f :\mathbb{R}^{n} \to \mathbb{R}\cup \{+\infty\}$. The regular (or Fr\'echet) subdifferential, limiting subdifferential and horizon subdifferential of $f$ at a point $\bar{\bm{x}} \in \textrm{dom}f$ respectively defined by
	\begin{equation*}
		\hat{\partial}f(\bar{\bm{x}}) := \left\{\bm{\xi} \in \mathbb{R}^{n}\mid  \liminf_{\substack{\bar{\bm{x}}\to\bm{x}\\ \bar{\bm{x}}\neq\bm{x}}} \frac{f(\bm{x}) - f(\bar{\bm{x}}) - \langle \bm{\xi}, \bm{x} -\bar{\bm{x}}\rangle}{\Vert \bm{x} -\bar{\bm{x}} \Vert}\geq  0\right\},
	\end{equation*}
	and 
	\begin{equation*}
		\partial f(\bar{\bm{x}}) := \left\{\bm{\xi} \in \mathbb{R}^{n}\mid \exists \bm{x}^{\nu}\xrightarrow{f}\bar{\bm{x}} \textrm{ and } \bm{\xi}^{\nu} \in \hat{\partial}f(\bm{x}^{\nu}) \textrm{ with } \bm{\xi}^{\nu}\to\bm{\xi}\right\},
	\end{equation*}
	and  
	\begin{equation*}
		\partial^{\infty} f(\bar{\bm{x}}) := \left\{\bm{\xi} \in \mathbb{R}^{n}\mid \exists \bm{x}^{\nu}\xrightarrow{f}\bar{\bm{x}} \textrm{ and } \bm{\xi}^{\nu} \in \hat{\partial}f(\bm{x}^{\nu}) \textrm{ with } \lambda^{\nu}\bm{\xi}^{\nu}\to\bm{\xi} \textrm{ for some } \lambda^{\nu} \searrow 0\right\} ,
	\end{equation*}
	where $\bm{x}^{\nu}\xrightarrow{f}\bar{\bm{x}}$ means $\bm{x}^{\nu} \to \bar{\bm{x}}$ and $f(\bm{x}^{\nu}) \to f(\bar{\bm{x}})$ and $\lambda^{\nu} \searrow 0$ means a sequence of real numbers $\lambda^{\nu} >0$ with $\lambda^{\nu} \to 0$. For any $\bm{\xi} \in \partial f(\bar{\bm{x}})$, we call it a \textit{subgradient} of $f$ at $\bar{\bm{x}}$.
\end{definition}
It follows from \cite[Theorem 8.6]{rockafellar2009variational} that $\hat{\partial}  f(\bar{\bm{x}})$, the set of regular subgradients of $f$ at $\bar{\bm{x}}$,
is closed and convex (though possibly empty). The set of subgradients, $\partial f(\bar{\bm{x}})$, is not
necessarily convex. If $f$ is a convex function, $\hat{\partial}   f(\bar{\bm{x}}) := \{ \bm{\xi} \in\mathbb{R}^{n} \mid f(\bm{x}) \ge f(\bar{\bm{x}}) + \la \bm{\xi}, \bm{x} - \bar{\bm{x}}\ra, \forall \ \bm{x} \in \mathbb{R}^n  \} = \partial f(\bar{\bm{x}})$ and reduces to the ordinary subdifferential of convex analysis \cite[Proposition 8.12]{rockafellar2009variational}. In particular, the subgradient of $\vert x\vert$ with respect to $x \in \mathbb{R}$ is $\xi = \text{sgn}(x)$ if $x\neq 0$; and $\xi \in [-1,1]$, otherwise. 

\begin{definition}
	Let $\mathcal{C}$ be a subset of $\mathbb{R}^{n}$ and let $\bar{\bm{x}} \in \mathcal{C}$ be given. The regular normal cone to $\mathcal{C}$ at $\bar{\bm{x}}$ is defined by
	\begin{equation*}\label{Def: Normal_cone}
		\widehat{N}_{\mathcal{C}}(\bar{\bm{x}}) := \left\{\bm{v}\in \mathbb{R}^{n}\mid \limsup_{\substack{\bm{x}\xrightarrow{\mathcal{C}}\bar{\bm{x}}\\\bm{x}\neq\bar{\bm{x}}}} \frac{\la\bm{v}, \bm{x}-\bar{\bm{x}}\ra}{\Vert\bm{x}-\bar{\bm{x}}\Vert_{2}} \leq 0\right\}.
	\end{equation*}
	Moreover, if there are sequences $\bm{x}^{\nu}\xrightarrow{\mathcal{C}}\bar{\bm{x}}$ and $\bm{v}^{\nu} \rightarrow \bm{v}$ with $\bm{v}^{\nu}\in\widehat{N}_{\mathcal{C}}(\bm{x}^{\nu})$, then the (general) normal cone to $\mathcal{C}$ at $\bar{\bm{x}}$ is written as $N_{\mathcal{C}}(\bar{\bm{x}})$.
\end{definition}

Under the assumption that $f:\mathbb{R}^{n}\to \mathbb{R}\cup\{+\infty\}$ is continuous and has at least one regular subgradient at $\bm{x}$, we know from \cite[Theorem 8.6]{rockafellar2009variational} that $\hat{\partial} f(\bm{x})$ is nonempty, closed, and convex. We embrace  the definition introduced in \cite[Theorem 3.6]{rockafellar2009variational}, and the horizon cone of $\hat{\partial} f(\bm{x})$ is defined by
\begin{equation}\label{eq:HorizonCone}
	\hat{\partial} f(\bm{x})^{\infty} := \{\bm{\xi}\mid \tilde{\bm{\xi}} + t \bm{\xi} \in \hat{\partial}f(\bm{x}),\ \forall \tilde{\bm{\xi}} \in \hat{\partial}f(\bm{x}), \forall t \in \mathbb{R}_{+}\}.
\end{equation}

\noindent Directly from the above definitions, we obtain from \cite[Theorem 8.6]{rockafellar2009variational} that
\begin{equation} \label{eq:relationship_Subdifferential}
	\hat{\partial}\Phi(\bar{\bm{x}}) \subseteq \partial\Phi(\bar{\bm{x}})\quad \textrm{ and } \quad 0 \in \hat{\partial}\Phi(\bar{\bm{x}})^{\infty} \subseteq \partial^{\infty}\Phi(\bar{\bm{x}}).
\end{equation}

Regularity is a key concept in nonsmooth analysis \cite{clarke1973necessary}. We recall the subgradient criterion for subdifferential regularity from \cite[Corollary 8.11]{rockafellar2009variational} and Clarke regularity of sets from \cite[Definition 6.4]{rockafellar2009variational} as follows.
\begin{definition}\label{Def_subdifferentiallyregular}
	Let $f:\mathbb{R}^{n} \to \mathbb{R}\cup\{+\infty\}$ and $\partial f(\bar{\bm{x}}) \neq \emptyset$ at a point $\bar{\bm{x}} \in \textrm{dom}f$. We say that $f$ is \textit{subdifferentially regular} at $\bar{\bm{x}}$ if 
	$$ \hat{\partial}f(\bar{\bm{x}}) = \partial f(\bar{\bm{x}}) \quad \textrm{ and } \quad  \hat{\partial}f(\bar{\bm{x}})^{\infty} = \partial^{\infty}f(\bar{\bm{x}}).$$
\end{definition}
\begin{definition}\label{Def_ClarkeRegularity}
	A set $\mathcal{C} \subseteq \mathbb{R}^{n}$ is \textit{regular} at one of its points $\bar{\bm{x}}$ in the sense of Clarke if it is locally closed at $\bar{\bm{x}}$ and every normal vector to $\mathcal{C}$ at $\bar{\bm{x}}$ is a regular normal vector, i.e., $N_{\mathcal{C}}(\bar{\bm{x}}) = \hat{N}_{\mathcal{C}}(\bar{\bm{x}})$.
\end{definition}

We next compute the regular, general, and horizon subgradients of $\Phi$ as follows.
\begin{proposition}\label{Prop:subdifferential}
	Consider \eqref{eq: main_Opt_General}. It holds for any $\bar{\bm{x}}$ that
	\begin{equation*}
		\begin{aligned}
			&\partial \Phi(\bar{\bm{x}})= \hat{\partial}\Phi(\bar{\bm{x}}) = \left\lbrace\bm{\xi}\mid \xi_i  \begin{cases}
				=\phi'(\vert\bar{x}_i\vert),& i \in \mathcal{I}(\bar{\bm{x}}),\\
				\in \begin{cases}
					\left[-\phi'(0),\phi'(0)\right],& \lim\limits_{t \to 0^{+}} \phi'(t) < +\infty,\\
					\mathbb{R},& \lim\limits_{t \to 0^{+}} \phi'(t) = +\infty,
				\end{cases}
				& i \in \mathcal{A}(\bar{\bm{x}}).
			\end{cases}\right\rbrace,\\
			&\hat{\partial} \Phi(\bar{\bm{x}})^{\infty}=\partial^{\infty}\Phi(\bar{\bm{x}})=\left\lbrace\bm{\xi}\mid \xi_i  \begin{cases}
				= 0,& i \in \mathcal{I}(\bar{\bm{x}}),\\
				\in \begin{cases}
					\{0\},& \lim\limits_{t \to 0^{+}} \phi'(t) < +\infty,\\
					\mathbb{R},& \lim\limits_{t \to 0^{+}} \phi'(t) = +\infty,
				\end{cases}
				& i \in \mathcal{A}(\bar{\bm{x}}).
			\end{cases}\right\rbrace.
		\end{aligned}
	\end{equation*}
	Consequently, the subdifferential regularity of $\Phi$ holds at every $\bm{x} \in \mathbb{R}^{n}$.
\end{proposition}
\begin{proof}
	Consider first $\phi$ on $\mathbb{R}$. Note that $\phi'(t) >0$ for any $t \in \mathbb{R}_{++}$ and $\phi'$ is subdifferentiable at $0$  under Assumption \ref{BasicAssum}. For any $i\in \mathcal{I}(\bar{\bm{x}})$, we have $\partial\phi(\vert \bar{x}_{i} \vert) = \hat{\partial}\phi(\vert \bar{x}_{i} \vert)=\phi'(\vert \bar{x}_{i} \vert)$ by \cite[Exercise 8.8]{rockafellar2009variational}. As for $i\in \mathcal{A}(\bar{\bm{x}})$, we consider the following two cases:
	\begin{itemize}
		\item[(i)] Consider $\lim\limits_{t \to 0^{+}} \phi'(t) < +\infty$. It follows from the definition of the regular subdifferential presented in Definition \ref{Def_Subdifferential} that $\hat{\partial}\phi(0) = \left[-\phi'(0),\phi'(0)\right]$. By \cite[Eq. 8(5)]{rockafellar2009variational}, one has $$\partial \phi(0) = \limsup\limits_{x_i\xrightarrow{\phi}0}\hat{\partial}\phi(\vert \bar{x}_i\vert) = \lim\limits_{x_i\xrightarrow{\phi}0}\phi'(\vert \bar{x}_i\vert) = \hat{\partial}\phi(0).$$
		\item[(ii)] Consider $\lim\limits_{t \to 0^{+}} \phi'(t) = +\infty$. This immediately gives that $\lim\limits_{x_i \to 0} \frac{\phi(\vert x_i \vert) - \phi(0)}{\vert x_i - 0 \vert} = +\infty$, hence indicating $$\liminf_{x_i \to 0} \frac{\phi(\vert x_i \vert) - \phi(0) - \xi_i (x_i-0)}{\vert x_i-0\vert} \geq 0$$ for any $\xi_{i} \in \mathbb{R}$, where the second inclusion holds by \cite[Theorem 8.6]{rockafellar2009variational}. Therefore, we have that $\mathbb{R} \subseteq \hat{\partial}\phi(0) \subseteq \partial\phi(0) \subseteq \mathbb{R}$. Hence, $\hat{\partial}\phi(0) = \partial\phi(0) =\mathbb{R}$.  
	\end{itemize}
	It then follows from \cite[Proposition 10.5]{rockafellar2009variational} that 
	$\hat{\partial}\Phi(\bar{\bm{x}}) = \partial\Phi(\bar{\bm{x}}) = \partial\phi(\vert \bar{x}_1 \vert) \times \ldots \times \partial\phi(\vert \bar{x}_n \vert)$, deriving  the desired results for $\partial\Phi(\bar\bx)$ and $\hat\partial\Phi(\bar \bx)$.

	As for $\hat\partial\Phi(\bar \bx)^\infty$,  since $\hat{\partial}\Phi(\bar{\bm{x}})$ is nonempty, closed and convex \cite[Theorem 8.6]{rockafellar2009variational}, it follows from \eqref{eq:HorizonCone} that 
	$\hat{\partial} \Phi(\bar{\bm{x}})^{\infty}$ takes the form as presented. 
	
	We next prove  $\hat{\partial}^{\infty}\Phi(\bar{\bm{x}})$ takes the same form as $\hat{\partial} \Phi(\bar{\bm{x}})^{\infty}$. By \eqref{eq:relationship_Subdifferential}, it suffices to show that $ \partial^{\infty} \Phi(\bar{\bm{x}}) \subset \hat{\partial}\Phi(\bar{\bm{x}})^{\infty}$. Let $\tilde{\bm{\xi}} \in \hat{\partial}\Phi(\bar{\bm{x}})$ and $t \geq 0$ be arbitrary. For any $\bm{\xi} \in \partial^{\infty} \Phi(\bar{\bm{x}})$, we know that there exist sequences $\{\bm{x}^{\nu} \}\xrightarrow{f}\bar{\bm{x}}$, $\{\lambda^{\nu}\}\searrow 0$, and $\bm{\xi}^{\nu}   \in \hat{\partial}\Phi(\bm{x}^{\nu} )$ such that $\{\lambda^{\nu}\bm{\xi}^{\nu} \} \to \bm{\xi}$.  Therefore,  $ \xi_i = \bm{0}$ since $\xi_i^\nu = \phi'(\vert \bar{x}_{i}\vert)$ is finite by Assumption \ref{BasicAssum}(ii) for $i \in \mathcal{I}(\bar{\bm{x}})$.  Hence, $\tilde{\xi}_i + t \xi_i = \tilde\xi_i$.  Similar arguments can also be applied to the cases where $i \in \mathcal{A}(\bar{\bm{x}})$, and we can deduce that 
	$\tilde{\xi}_i + t \xi_i = \tilde\xi_i$ for $\lim_{t\to 0^+}\phi'(t) < +\infty$ 
	and $\tilde{\xi}_i + t \xi_i \in\mathbb{R}$ for $\lim_{t\to 0^+}\phi'(t) = +\infty$. Overall,   
	$\tilde{\bm{\xi}} + t \bm{\xi} \in \hat{\partial}\Phi(\bar{\bm{x}})$.
	
	By Definition \ref{Def_subdifferentiallyregular}, we know $\Phi$ is subdifferentially regular at any $\bm{x} \in \mathbb{R}^{n}$.
	 \end{proof}

The regular and general normal vectors associated with $\Omega$ in \eqref{eq: main_Opt_General} are calculated as follows.

\begin{proposition}\label{Prop_NormalCone}
	Consider \eqref{eq: main_Opt_General}. It holds for any $\bar{\bm{x}} \in \textbf{bd}\ \Omega$ that
	\begin{itemize}
		\item[(i)] $\Omega$ is regular in the sense of Clarke at $\bar{\bm{x}}$, i.e., $N_{\Omega}(\bar{\bm{x}}) = \widehat{N}_{\Omega}(\bar{\bm{x}})$.
		\item[(ii)] $N_{\Omega}(\bar{\bm{x}}) = \textrm{pos}\ \partial\Phi(\bar{\bm{x}}) \cup \partial^{\infty}\Phi(\bar{\bm{x}})$. Indeed, $N_{\Omega}(\bar{\bm{x}}) = \lambda \partial \Phi(\bar \bx) \ \textrm{for any}\ \lambda \ge 0$. 
	\end{itemize}
\end{proposition}
\begin{proof}
	For statement (i), it follows from Definition \ref{Def_subdifferentiallyregular} and Proposition \ref{Prop:subdifferential} that $\Phi$ is regular at $\bar{\bm{x}}$ with $0 \notin \partial \Phi(\bar{\bm{x}})$. Then $\Omega$ is regular at $\bar{\bm{x}}$ by \cite[Proposition 10.3]{rockafellar2009variational}. Hence, by Definition \ref{Def_ClarkeRegularity}, we know that $N_{\Omega}(\bar{\bm{x}}) = \widehat{N}_{\Omega}(\bar{\bm{x}})$. Statement (ii) can be obtained straightforwardly from the definition of positive hull and Proposition \ref{Prop:subdifferential} by \cite[Proposition 10.3]{rockafellar2009variational}.
	 \end{proof}
The following corollary explicitly presents the computation of the normal cone to the constraint $\mathcal{B}_{\ell_p}$, whose proof can also be directly referred to \cite[Theorem 2.2]{wang2021constrained}.
\begin{corollary}\label{eq: Concrete_elements}
	Consider \eqref{eq: main_Opt}. For any $\bm{x} \in \mathcal{B}_{\ell_p}$, the normal cone to $\mathcal{B}_{\ell_p}$ at $\bm{x}$ has the form:
	\begin{equation*}
		N_{\mathcal{B}_{\ell_p}}(\bm{x}) = \left\lbrace 
		\begin{aligned}
			&\{\bm{v}\in \mathbb{R}^{n}\mid v_i = \lambda\text{sgn}(x_i)p\vert x_i\vert^{p-1}, i\in\I(\bm{x});\lambda \geq 0\}, &\textrm{ if } \bm{x} \in \textbf{bd}\ \mathcal{B}_{\ell_p}, \\
			&\{\bm{0}\}, &\textrm{ if } \bm{x} \in \textbf{int}\ \mathcal{B}_{\ell_p}.
		\end{aligned}\right. 
	\end{equation*}
\end{corollary}


To characterize optimality conditions for \eqref{eq: main_Opt_General}, we define the following constraint qualification on the constraint set $\Omega$.

\begin{definition}[Nonsmooth MFCQ for \eqref{eq: main_Opt_General} {\cite[Eq. (3.3)]{dempe2011generalized}}]\label{Def_EMFCQ}
	Consider \eqref{eq: main_Opt_General} under Assumption \ref{BasicAssum}. For a feasible point $\bm{x}$ of \eqref{eq: main_Opt_General}, we say that the nonsmooth Mangasarian-Fromovitz constraint qualification (MFCQ) is satisfied at $\bm{x} \in \Omega$ if the following implication holds
	$$\left[\bm{0} \in \partial(\lambda\vartheta)(\bm{x}), \lambda \in N_{\mathbb{R}_{-}}(\vartheta(\bm{x}))\right] \Longrightarrow \lambda = 0,$$  where $\vartheta(\bm{x}) = \Phi(\bm{x}) - \gamma$.
\end{definition}

The nonsmooth MFCQ holds naturally true in the interior of $\Omega$ since $\lambda = 0$ in this case. The following result shows it also holds true on the boundary of $\Omega$.

\begin{proposition}\label{pro_MFCQholds} 
	The nonsmooth MFCQ naturally holds at each point in $\Omega$ under Assumption \ref{BasicAssum}.
\end{proposition}
\begin{proof} 
	We only consider the case where $\bar x \in \textbf{bd}\ \Omega$. We have that 
	\begin{equation}\label{normal.lambda}
		\partial \delta_{\mathbb{R}_-}(\vartheta(\bar{\bm{x}})) = N_{\mathbb{R}_-}(\vartheta (\bar{\bm{x}}))
		= N_{\mathbb{R}_-}(0) =\{ \lambda \mid  \lambda \ge 0\}. 
	\end{equation}
	On the other hand, the nonnegative rescaling property of subgradients (see \cite[Eq. 10(6)]{rockafellar2009variational}), together with Proposition \ref{Prop:subdifferential}, gives
	\begin{equation}\label{eq: Generalresults}
		\begin{aligned}
			&\quad\ \partial(\lambda \vartheta)(\bar{\bm{x}}) 
			= \partial(\lambda \Phi)(\bar{\bm{x}})
			= \lambda \partial\Phi(\bar{\bm{x}}).
		\end{aligned}
	\end{equation}
	From Proposition \ref{Prop:subdifferential}, 
	$\lambda \partial\Phi(\bar{\bm{x}}) \ne \bm{0} $ unless 
	$\mathcal{I}(\bar{\bm{x}}) = \emptyset$, which cannot happen since 
	$\bar x \in \textbf{bd}\ \Omega$. 
	Therefore, $\lambda = 0$, meaning the nonsmooth MFCQ holds true in this case. 
	 \end{proof}

Now we can characterize the first-order stationary conditions to problem \eqref{eq: main_Opt_General}. 
\begin{proposition}\label{GelLocalImpliesSP}
	Let $\bm{x}^{*} \in \textbf{bd}\ \Omega$ be a local minimizer of problem \eqref{eq: main_Opt_General}. Then there exists $\lambda^{*} \in \mathbb{R}_{+}$ such that the following conditions are satisfied:
	\begin{equation}\label{Def SC} 
		\begin{cases}
			i \in \mathcal{I}(\bm{x}^{*}):  \nabla_i f(\bm{x}^{*}) + \lambda^{*} \phi'(\vert x_{i}^{*}\vert)\text{sgn}(x_{i}^{*}) = 0,\\
			i\in \mathcal{A}(\bm{x}^{*}): 
			\begin{cases}
				\vert \nabla_i f(\bm{x}^{*})\vert \leq \lambda^{*}\phi'(0),  & \lim\limits_{t \to 0^{+}} \phi'(t) < +\infty,\\
				-\nabla_i f(\bm{x}^{*}) \in \mathbb{R}, & \lim\limits_{t \to 0^{+}} \phi'(t) = +\infty;
			\end{cases}\\
			\Phi(\bm{x}^{*}) - \gamma = 0, \ \lambda^{*} \geq 0.		
		\end{cases}  
	\end{equation}
	\end{proposition}
	\begin{proof}
It follows from \cite[Exercise 10.52]{rockafellar2009variational} (with $X = \mathbb{R}^{n}$ and $\theta =  \delta_{\mathbb{R}_{-}}$) that there exists a scalar $\lambda^{*}$ such that
\[
\bm{0} \in \partial (f + \lambda^{*}\vartheta)(\bm{x}^{*}),\quad 
\lambda^{*} \in \partial \delta_{\mathbb{R}_{-}}(\vartheta(\bm{x}^{*})). \]
Since $\bm{x}^{*} \in \textbf{bd}\ \Omega$, we obtain from \eqref{normal.lambda} that 
\begin{equation}\label{eq: Var_condition}
	\begin{aligned}
		&\quad\ \bm{0} \in \partial (f + \lambda^{*}\vartheta)(\bm{x}^{*}) = \nabla f(\bm{x}^{*}) + \partial (\lambda^{*} \vartheta)(\bm{x}^{*}) = 
		\nabla f(\bm{x}^{*}) +  \lambda^{*} \partial  \vartheta(\bm{x}^{*}),
	\end{aligned}
\end{equation}
where the first equality follows from \cite[Exercise 8.8(c)]{rockafellar2009variational} and the second equality holds by \eqref{eq: Generalresults}. Then
combining \eqref{eq: Var_condition} and Proposition \ref{Prop:subdifferential} completes the 
proof of \eqref{Def SC}. 
\end{proof}
\begin{corollary}\label{LocalImpliesSP}
Let $\bm{x}^{*} \in \textbf{bd}\ \mathcal{B}_{\ell_{p}}$ be a local minimizer of problem \eqref{eq: main_Opt}. Then there exists $\lambda^{*} \in \mathbb{R}_{+}$ such that one of the following conditions is satisfied:
\begin{equation}\label{Def SC_lp} 
	\begin{cases}
		\nabla_i f(\bm{x}^{*}) + \lambda^{*} p\vert x_i^{*}\vert^{p-1}\text{sgn}(x_{i}^{*}) = 0, &\ i \in \I(\bm{x}^{*}); \\
		\Vert\bm{x}^{*} \Vert_{p}^p - \gamma = 0, \ \lambda^{*} \geq 0.	
	\end{cases}  
\end{equation}
\end{corollary}
\noindent In virtue of Corollary \ref{LocalImpliesSP}, we therefore say that an $\bm{x}^{*} \in \mathcal{B}_{\ell_p}$ is a \textit{stationary point} of \eqref{eq: main_Opt} if there exists $\lambda^{*} \in \mathbb{R}_{+}$ such that $(\bm{x}^{*},\lambda^{*})$ satisfies \eqref{Def SC_lp} or $\nabla f(\bm{x}^*) = \bm{0}$.

\section{A Hybrid First-order Algorithm for Solving \eqref{eq: main_Opt}}\label{Algo_sec}
In this section we present  a hybrid first-order algorithm for solving \eqref{eq: main_Opt}. This algorithm adaptively alternates between solving a Frank-Wolfe (FW) subproblem and a gradient projection (GP) subproblem, depending upon the relative relationship between the current solution estimate and the constraint set. As highlighted in the introduction, each iteration of the proposed algorithm involves only one type of subproblems. And for both types of subproblems many efficient solvers can be readily leveraged  within the proposed algorithmic framework.

\subsection{Main Algorithm Framework}
We formally present our proposed algorithm in Algorithm \ref{alg.propose}, consisting of two computational building blocks. The first block encompasses the Frank-Wolfe subproblem with backtracking line search, while the second addresses convex subproblems of the Euclidean projections onto weighted $\ell_{1}$ balls. 

\begin{algorithm}[H]
	\caption{A Hybrid First-order Method for Solving \eqref{eq: main_Opt}}
	\label{alg.propose}
	\begin{algorithmic}[1]
		\State \textbf{Initialization:} $\bm{x}^{0} \in \mathcal{B}_{\ell_{p}}$, $p\in(0,1)$, $\gamma > 0$, $\delta_{\textrm{FW}} \geq 0$, $\delta_{\textrm{PG}} \geq 0$, $\alpha_{\max} = 1$, $\beta \in (0,1/L_f)$ and initial Lipschitz estimate $L_{-1} > 0$
		\For{$k = 0, 1, ...$ }
		\If{$\bm{x}^{k} \in \textbf{int}\ \mathcal{B}_{\ell_{p}}$} \label{line:InteriorCase}
		\State  Compute $\bm{s}^{k} \in \mathcal{B}_{\ell_{p}}$ from \eqref{eq: FW_linearization} \label{line:interMediatepoint} 
		\State Compute $\bm{d}^{k} = \bm{s}^{k} - \bm{x}^{k}$ \Comment{(\small the FW direction)} \label{line:updateDirection}
		\State Compute $g_{k} := \la-\nabla f(\bx^{k}), \bm{d}^{k}\ra$ \Comment{(\small FW gap by Lemma \ref{Lem: subproblem_property}(i))} \label{line:GapValue}
		\If{ $g_{k} \leq \delta_{\text{FW}}$}
		\textbf{return} $\bm{x}^{k}$ \label{line:OptFW}
		\EndIf
		\State $\alpha^{k}, L_{k} = \texttt{STEP\_SIZE}(f,\bm{d}^{k},\bm{x}^{k},g_{k},L_{k-1},\alpha_{\max})$ \label{line:StepSize}
		\State $\bm{x}^{k} \gets \bx^{k} + \alpha^{k}\bm{d}^{k}$.
		\ElsIf{$\bm{x}^{k}\in \textbf{bd}\ \mathcal{B}_{\ell_{p}}$} \label{line:bounadryCase}
		\State Compute $\bm{u}^{k} = \bm{x}^{k} - \beta\nabla f(\bm{x}^{k})$ \label{line:gradStep}\Comment{(\small the gradient step)} 
		\State Compute $\mathcal{I}^{k}$ and $\mathcal{A}^{k}$ \label{line:UpdateIndex}
		\State Compute $\gamma_{\text{PG}}^{k} = p\gamma$ \label{line:radiusOFpgd}
		\State Compute $\bm{x}^{k+1}_{\I^{k}} = \texttt{PROJ}(\bm{x}^{k}_{\I^{k}},\bm{u}^{k}_{\I^{k}},\gamma_{\text{PG}}^k)$ and $\bm{x}^{k+1}_{\mathcal{A}^{k}} = \bm{0}_{\mathcal{A}^{k}}$ from \eqref{sub.proj} \qquad\Comment{\text{(\small the projection step)}}\label{line:projStep}
		\State Compute $r^{k} := \Vert \bm{x}^{k+1} - \bm{x}^{k} \Vert_2^2$  \Comment{(\small PG gap by Lemma \ref{Lem: subproblem_property}(ii))}
		\If{$r^{k}\leq \delta_{\text{PG}}$}
		\textbf{return} $\bm{x}^{k+1}$ \label{line:OptPG}
		\EndIf
		\State $\bm{x}^{k} \gets \bm{x}^{k+1}$
		\EndIf
		\EndFor
	\end{algorithmic}
\end{algorithm}
\subsubsection{Frank-Wolfe Block}\label{FW_subsec}
In this case, the condition in Line \ref{line:InteriorCase} holds. Our basic idea is to approximate the objective function $f$ through a local linear expansion while preserving the nonconvex $\ell_{p}$ ball constraint in \eqref{eq: main_Opt}. Therefore, at current iterate $\bm{x}^{k}$, the following Frank-Wolfe type subproblem is solved for an intermediate point $\bm{s}^{k} \in \mathcal{B}_{\ell_{p}}$ in  Line \ref{line:interMediatepoint}:  
\begin{equation}\label{eq: FW_linearization}\tag{\text{$\mathscr{P}_{\textrm{FW}}$}}
	\begin{aligned}
		\min_{\bm{s}\in\mathbb{R}^{n}}\ Q(\bm{s};\bm{x}^{k}) := \langle \nabla f(\bm{x}^{k}), \bm{s}\rangle \quad
		\text{s.t.} \ \bm{s} \in \mathcal{B}_{\ell_{p}}. 
	\end{aligned} 
\end{equation}
Although the feasible set in \eqref{eq: FW_linearization} makes it nonconvex,  it indeed admits a closed-form solution as shown in the following theorem, and its proof can be found in Appendix \ref{app_Sub}.
\begin{theorem}\label{Theo: FW_reformulation} 
	Suppose $\nabla f(\bm{x}^{k}) \neq \bm{0}$. Let $\I^k_{\max}  := \argmax\limits_{i\in[n]}\vert\nabla_{i} f(\bm{x}^{k})\vert$. A global optimal solution $\bm{s}^{k}$ of \eqref{eq: FW_linearization} is given by:
	\begin{equation}\label{eq: FW_direction}
		s_i^{k}=  
		\begin{cases}
			-\textrm{\text{sgn}}(\nabla_{i} f(\bm{x}^{k}))\gamma^{\frac{1}{p}}, &\text{if  $i\in \I^k_{\max}$},  \\
			0, &\text{otherwise}.
		\end{cases}\\
	\end{equation}
\end{theorem}
\begin{remark}
	The proof of Theorem \ref{Theo: FW_reformulation} is intuitively simple. 
	Note that the objective of \eqref{eq: FW_linearization} is linear with respect to $\bm{s}$, it is easy to see that the global optimal solution should be located at a vertex of a $\ell_{p}$ ball. Based on this key observation, it is clear that $\bm{s}^{k} \in \argmin_{\bm{s} \in \mathcal{B}_{\ell_p}} \langle \nabla f(\bm{x}^{k}),\bm{s}\rangle = \argmin_{\bm{s} \in \textrm{conv}(\mathcal{B}_{\ell_p})} \langle \nabla f(\bm{x}^{k}),\bm{s}\rangle$, in which $\textrm{conv}(\mathcal{B}_{\ell_p})$ is exactly an $\ell_{1}$ ball. Therefore, we can directly extract the solution from a counterpart $\ell_1$ ball-constrained problem with an appropriate scaling of the radius $\gamma$.
\end{remark}

\noindent With $d^k$, a FW direction $\bm{d}^{k} = \bm{s}^{k} - \bm{x}^{k}$ is computed in Line \ref{line:updateDirection} and the FW gap, which can be used to represent the optimality measure by Lemma \ref{Lem: subproblem_property}(i), is computed in Line \ref{line:GapValue}. If the condition $g^{(k)} \leq \delta_{\text{FW}}$ in Line \ref{line:OptFW} is not satisfied,  a backtracking line-search along  $\bm{d}^{k}$, together with a bisection method (subroutine $\texttt{STEP\_SIZE}$) is performed in Line \ref{line:StepSize}, as described in  \S\ref{subsection_FW_stepsize}.

\subsubsection{Gradient Projection Block}\label{PGM_subsec_}
In this case, the condition in Line \ref{line:bounadryCase} holds and subsequent steps perform projections onto a weighted $\ell_{1}$ ball in a reduced space. We approximate the objective $f$ and the feasible set $\mathcal{B}_{\ell_{p}}$ in \eqref{eq: main_Opt} by their quadratic majorant and affine majorant respectively at the current estimate $\bm{x}^{k}$. Specifically, a direct consequence of Assumption \ref{BasicAssum}(i) is that $f$ admits a quadratic upper bound for $L_f > 0$ and any feasible iterates, i.e.,
\begin{equation}\label{eq:QuadrticBound}
	f(\bm{x}) \leq f(\bm{x}^{k}) + \la \nabla f(\bm{x}^{k}),\bm{x} - \bm{x}^{k}\ra + \frac{L_{f}}{2}\Vert \bm{x} - \bm{x}^{k}\Vert, \ \forall \bm{x} \in \mathcal{B}_{\ell_{p}}.
\end{equation}
Hence, we turn to minimizing the surrogate function, given in the right-hand side of \eqref{eq:QuadrticBound}. On the other hand, by the concavity of the $\ell_{p}$ norm, we have 
\begin{equation}\label{eq: well-posedness}
	\Vert\bx\Vert_p^p \le \Vert\bx^k\Vert_p^p + \sum_{i\in \I^k} w_{i}^{k}(\vert x_{i}\vert - \vert x_{i}^{k}\vert) \le \gamma = \Vert\bx^k\Vert_p^p, \ \forall \bm{x} \in \mathcal{B}_{\ell_{p}},
\end{equation}
where $w_{i}^{k} = p\vert x_{i}^{k}\vert^{p-1},\ i\in \I^{k}$. It follows from \eqref{eq: well-posedness} that $\sum_{i\in \I^k} w_{i}^{k}(\vert x_{i}\vert - \vert x_{i}^{k}\vert) \leq 0$. We then replace  $\mathcal{B}_{\ell_{p}}$ by the  weighted $\ell_{1}$ ball
$$ \mathcal B_{\ell_1}^k := \{\bm{x} \in \mathbb{R}^{n} \mid \langle \bm{w}_{\I^k}^{k} , \vert \bm{x}_{\I^k} \vert - \vert \bm{x}^{k}_{\I^k} \vert \rangle  \leq 0, \textrm{sgn}(\bm{x}^{k}_{\mathcal{I}^k}) \circ \bm{x}_{\mathcal{I}^k} \geq \bm{0}, \bm{x}_{\mathcal{A}^{k}} = \bm{0}\}.$$ 

In Line \ref{line:gradStep} we compute the point  $\bm{u}^{k} = \bm{x}^{k} - \beta\nabla f(\bm{x}^{k})$  that will be projected onto a weighted $\ell_{1}$ ball whose  radius is calculated in Line \ref{line:radiusOFpgd}. Therefore, in Line \ref{line:projStep}, the following projection subproblem is solved for $\bx^{k+1}$  if $\bx^k$ is not optimal to \eqref{eq: main_Opt}
\begin{equation}\label{sub.proj}\tag{\text{$\mathscr{P}_{\textrm{GP}}$}}
	\begin{aligned}
		\min_{\bm{x} \in 
			\mathbb{R}^{n}} \ &\quad P(\bm{x};\bm{x}^{k}):=  \frac{1}{2\beta} \Vert\bm{x} - \bm{u}^{k}\Vert_2^{2}\\
		\text{s.t.}\ & \quad {\bm x}\in \mathcal B_{\ell_1}^k,
	\end{aligned}
\end{equation}
where the stepsize $\beta > 0$. In particular, if $\beta = 1 / L_f$,  the objective $P$  is the right-hand side of \eqref{eq:QuadrticBound}. Here, we simply use a constant stepsize $\beta \in (0, 1/L_f)$ which by \cite{goldstein1964convex} can guarantee sufficient decrease in the objective caused by $\bx^{k+1}$. 
The new iterate stays in the same orthant of $\bx^k$ 
and keeps the zero components in $\bx^k$ still as zero by constraints  $\bm{x}_{\mathcal{A}^{k}} = \bm{0}$ and $\textrm{sgn}(\bm{x}^{k}_{\mathcal{I}^k}) \circ \bm{x}_{\mathcal{I}^k} \geq \bm{0}$. This requirement ensures that the iterates have the same support when approaching the optimal solution. The weighted projection subproblem \eqref{sub.proj} is solved by subroutine $\texttt{PROJ}$.

\begin{remark}[Well-posedness of \eqref{sub.proj}]
	It should be noted that the subproblem \eqref{sub.proj} is well-defined in the sense that the feasible set $\mathcal{B}_{\ell_1^{\bm{w}^k}}$ is nonempty and is a subset of $\mathcal{B}_{\ell_{p}}$. To see this, we have
	\begin{equation*}\label{eq:well-posedness}
		\begin{aligned}
			\bm{x} \in  \mathcal{B}_{\ell_1^{\bm{w}^k}} &= \{\bm{x} \in \mathbb{R}^{n} \mid \langle \bm{w}_{\I^k}^{k} , \vert \bm{x}_{\I^k} \vert - \vert \bm{x}^{k}_{\I^k} \vert \rangle  \leq 0, \textrm{sgn}(\bm{x}^{k}_{\mathcal{I}^k}) \circ \bm{x}_{\mathcal{I}^k} \geq \bm{0}, \bm{x}_{\mathcal{A}^{k}} = \bm{0}\}\\
			&\subseteq \{\bm{x} \in \mathbb{R}^{n} \mid \sum_{i=1}^{n}\vert { x}_i \vert^{p} \leq \gamma\} = \mathcal{B}_{\ell_{p}}  \textrm{ with } \gamma > 0,
		\end{aligned}
	\end{equation*}
	where the set inclusion holds by \eqref{eq: well-posedness}.  Indeed, $\bm{x}^{k} \in \textbf{bd}\ \mathcal{B}_{\ell_p}$ is a Slater point of the feasible set because $\bm{x}^{k} \in \mathcal{B}_{\ell_1^{\bm{w}^k}}$. Moreover, it follows from the compactness of $\mathcal{B}_{\ell_{p}}$ that the optimal solution set of \eqref{sub.proj} is nonempty. Under Slater's condition and the positivity and finiteness of $\bm{w}_{\mathcal{I}^k}$, \eqref{sub.proj} can be effectively and efficiently solved by many algorithms with complexity $O(n)$ in practice \cite{condat2016fast,perez2022efficient}. These arguments establish the well-posedness of \eqref{sub.proj}.
\end{remark}

\begin{remark}
	Algorithm \ref{alg.propose} can be appropriately extended to solve \eqref{eq: main_Opt_General} covering a broader class of nonconvex sparsity-promoting regularizers. For this purpose, one needs to make the computation of the FW step clear and specify the weights in the PG step, and all other computations follow the similar spirit of Algorithm \ref{alg.propose}. Thanks to the linear nature of the objective function and  the function $\phi^{-1}:\mathbb{R}_{+} \to \mathbb{R}_{+}$ is continuously
	differentiable, a global optimal solution $\bm{s}^{k}$ is given by:
	\begin{equation*}\label{eq: FW_direction_General}
		s_i^{k}=  
		\begin{cases}
			-\text{sgn}(\nabla_{i} f(\bm{x}^{k}))\phi^{-1}(\gamma), &\text{if } i\in \I^k_{\max},  \\
			0, &\text{otherwise}.
		\end{cases}\\
	\end{equation*}
	To generate the PG step, we simply need to set the weights $w_i^{k}, i \in \mathcal{I}(\bm{x}^{k})$ at the $k$th iteration as
	\begin{equation*}\label{eq: weights_cal}
		w_{i}^{k} = w(x_{i}^{k}) = \phi'(x_i),\ \forall i \in \mathcal{I}(\bm{x}^{k}).
	\end{equation*}
\end{remark}

\subsection{Computing the FW Subproblem Stepsize}\label{subsection_FW_stepsize} 
Line \eqref{line:StepSize} of Algorithm \ref{alg.propose} corresponds to a backtracking line-search procedure, with details presented in Algorithm \ref{alg.subprocedure_fw}. 
This procedure computes a stepsize estimate depending upon the local properties of the objective $f$, as opposed to having access to the knowledge of its global Lipschitz constant $L_f$. Specifically, the following quadratic approximation to $f({\bm x}^k + \alpha {\bm d^k})$ is adopted as the line-search objective, namely,

\begin{equation}\label{eq: QuadrticApproximation}
	f^{\textrm{surro}}_{k}(\alpha,M) := f(\bm{x}^{k}) - \alpha g_{k} + \frac{\alpha^2 M}{2}\Vert \bm{d}^{k}\Vert^{2}_{2}, 
\end{equation}
where $\alpha \in \left[ 0,\alpha_{\max}\right]$ and $M > 0$ serves as an (local) estimate of $L_f$ and 
gives an initial stepsize estimate, as shown in Line \ref{line:localiniStep}. In the $\texttt{While}$ loop (corresponding to Lines \ref{line:Begin_Backtracking}-\ref{line:End_Backtracking}), if the sufficient decrease condition in Line \ref{line:Begin_Backtracking} is not satisfied, then $M$ increases by a fraction of $\tau>1$, and the stepsize is reduced accordingly. 

Lines \ref{line:Begin_rootfinding}-\ref{line:End_rootfinding} correspond to the root-finding block, and this block is performed when such an $\alpha$ obtained in Line \ref{line:stepsizeByBT} leads us to move out of $\mathcal{B}_{\ell_{p}}$.   Therefore, if the condition in Line \ref{line:Begin_rootfinding} is triggered, we select an $\alpha_{\textrm{bis}}$ to ensure the new iterate $\bm{x}^{k+1} \in \textbf{bd}\ \mathcal{B}_{\ell_{p}}$, i.e., finding an $\alpha_{\textrm{bis}}$ which is a root of
\begin{equation}\label{eq: univariate_equation}
	\chi(\alpha) = \Vert \alpha\bm{s}^{k} + (1-\alpha)\bm{x}^{k}\Vert_{p}^{p} - \gamma = 0
\end{equation}
on $(0, \alpha)$   in Line \ref{line:stepsizeByBT}. Note that this univariate nonlinear equation always has a solution over $(0,\alpha)$ since $\chi(0) < 0$ and $\chi(\alpha) > 0$. The bisection method is employed to determine the root of the nonlinear equation \eqref{eq: univariate_equation}.

\begin{algorithm}[htbp]
	\caption{Backtracking Line-search and Bisection in the FW Block}
	\label{alg.subprocedure_fw}
	\begin{algorithmic}[1]
		\Procedure{$\texttt{step\_size}(f,\bm{d}^{k},\bm{x}^{k},g_{k},L_{k-1},\alpha_{\max})$}{}
		\State Choose $\tau > 1$, $\zeta \leq 1$ and $\alpha_{\textrm{bis}}>0$\label{line:LocalCons}
		\State Choose $M \in \left[\zeta L_{k-1},L_{k-1}\label{line:LocalPara} \right]$\label{line:localLip}
		\State $\bar{\alpha}^k = \min\{g_{k}/(M\Vert \bm{d}^{k}\Vert^{2}),\alpha_{\max}\}$\label{line:localiniStep}
		\While{$f(\bm{x}^{k}+\bar{\alpha}^k\bm{d}^{k}) > f^{\textrm{surro}}_{k}(\bar{\alpha}^k,M)$} \label{line:Begin_Backtracking}
		\State $M = \tau M$ \label{line:tauMultiplication}
		\State $\bar{\alpha}^k = \min\{g_{k}/(M\Vert \bm{d}^{k}\Vert^{2}),\alpha_{\max}\}$ \label{line:stepsizeByBT}
		\EndWhile\label{line:End_Backtracking}
		\State $\alpha \gets \bar{\alpha}^k$
		\If{$\chi(\alpha) > \delta_{\textrm{bis}}$} \label{line:Begin_rootfinding}
		\State Find an $\alpha_{\textrm{bis}}$ such that $\chi(\alpha_{\textrm{bis}}) \leq \delta_{\textrm{bis}}$ over $(0,\alpha)$
		\State $\alpha \gets \alpha_{\textrm{bis}}$
		\EndIf \label{line:End_rootfinding}
		\State \Return $\alpha$, $M$
		\EndProcedure
	\end{algorithmic}
\end{algorithm}

\begin{remark}
	The step-size $\alpha_{\textrm{bis}}$ results in $\bm{x}^{k+1} \in \textbf{bd}\ \mathcal{B}_{\ell_{p}}$, thereby triggering a gradient projection step. 
	The reason that we do not construct and solve a Frank-Wolfe subproblem at an iterate $\bm{x}^k \in \textbf{bd}\ \mathcal{B}_{\ell_{p}}$ is to avoid creating an infeasible search direction---a difficult situation for convex constrained problems. To see this, suppose $\bm{x}^k \in \textbf{bd}\ \mathcal{B}_{\ell_{p}}$ in  \eqref{eq: FW_linearization}, which generates a vertex point (by Theorem \ref{Theo: FW_reformulation}) $\bm{s}^k$ of 
	$\mathcal{B}_{\ell_{p}}$.  
	Then any point $\bm{x}(\alpha) = (1-\alpha) \bm{x}^k + \alpha \bm{s}^k $ with $\alpha\in(0,1)$   could be infeasible if $\bm{s}^k$ and $\bm{x}^k$ are in the same orthant, i.e., $\text{sgn}(x_i)\text{sgn}(s_i) \ge 0, \forall i$. In fact, for any   $i\in \{ i \mid  s^k_i \ne 0 \text{  or } x^k_i \ne 0\} $ (this set is nonempty since $\bm{s}^k, \bm{x}^k \in \textbf{bd}\ \mathcal{B}_{\ell_{p}}$),  
	$$\begin{aligned} x_i(\alpha) = \vert (1-\alpha) x^k_i + \alpha s^k_i \vert^p  
		=   \vert (1-\alpha) |x^k_i| + \alpha |s^k_i| \vert^p  
		>   (1-\alpha) \vert x^k_{i }\vert^p + \alpha \vert s^k_{i} \vert^p
	\end{aligned}$$
	for any $\alpha \in (0,1)$  by the strict concavity of $\vert(\cdot)\vert^{p}$.  It then follows that $$\|\bm{x}(\alpha)\|_p^p =   \Vert(1-\alpha)\bm{x}^{k} + \alpha\bm{s}^{k}\Vert_{p}^{p}
	>  (1 - \alpha)\Vert\bm{x}^{k}\Vert_{p}^{p} + \alpha\Vert\bm{s}^{k}\Vert_{p}^{p}
	= \gamma,$$ meaning $\bm{x}(\alpha)$ is infeasible for any $\alpha\in(0,1)$. 
\end{remark}

\section{Convergence Analysis}\label{Conver_sec} 
In this section, we establish the global convergence of Algorithm \ref{alg.propose} by setting $\delta_{\textrm{FW}}=\delta_{\textrm{GP}}=0$. Since Algorithm \ref{alg.propose} adaptively alternates between solving a FW subproblem and a projection subproblem,  we can split the sequence generated by the algorithm into two subsequences according to the type of the subproblems.  We then show the global convergence of each subsequence.    
As a result,  each cluster point of the subsequences satisfies the optimality conditions of \eqref{eq: main_Opt}.  
For this purpose, we first define  
\begin{equation*}
	\begin{aligned}
		\mathcal{S}_1 &= \{k \in \mathbb{N} \mid  \Vert\bx^k\Vert_p^p < \gamma \}    \quad  \text{and}\quad 
		\mathcal{S}_2 &= \{k \in \mathbb{N} \mid  \Vert\bx^k\Vert_p^p = \gamma \}.
	\end{aligned}
\end{equation*}
For $k\in \mathcal{S}_1$, we further define $\mathcal{S}_1 = \mathcal{T}_1 \cup \mathcal{T}_2$ with 
\begin{equation*}
	\begin{aligned}
		\mathcal{T}_1 = \{k\in \mathcal{S}_1\mid \alpha^k = \bar\alpha^k \}    \quad  \text{and}\quad 
		\mathcal{T}_2 = \{k\in \mathcal{S}_1 \mid \alpha^{k} =   \alpha^k_\textrm{bis}\}.
	\end{aligned}
\end{equation*}
For ease of presentation, at the $k$th iterate, define the reduction in $f$ and subproblem objective $P$ caused by $\bm{x}^{k+1}$
as 
$$\Delta f(\bm{x}^{k+1}) = f(\bm{x}^{k}) - f(\bm{x}^{k+1}),\quad \Delta P(\bm{x}^{k+1};\bm{x}^{k}) = P(\bm{x}^{k};\bm{x}^{k}) - P(\bm{x}^{k+1};\bm{x}^{k}),$$ 
and the reduction in $Q$ by $\bm{s}^{k}$ as
$$\Delta Q(\bm{s}^{k}; \bm{x}^{k}) = Q(\bm{x}^{k};\bm{x}^{k}) - Q(\bm{s}^{k};\bm{x}^{k}).$$

\subsection{Optimality Conditions of Subproblems}
Before proceeding, we first present the optimality conditions of \eqref{eq: FW_linearization} and \eqref{sub.proj} accordingly. Let $\bm{s}^{k}$ be the optimal solution of \eqref{eq: FW_linearization} at the $k$th iteration. Then according to Corollary \ref{LocalImpliesSP}, 
there exists $\varsigma^k \in \mathbb{R}_{+}$ such that 
\begin{subequations}\label{eq:optimality_FW_sub}
	\begin{alignat}{1}
		\nabla_if(\bm{x}^k) + \varsigma^k p \vert s_i^{k}\vert^{p-1}\textrm{sgn}(s_i^{k}) &= 0,\ i\in \mathcal{I}(\bm{s}^{k}),\\
		\sum_{i \in \mathcal{I}(\bm{s}^{k})} \vert s_i^{k}\vert^{p} &= \gamma.
	\end{alignat}
\end{subequations}
For the optimal solution, ${\bm x}^{k+1}$, of \eqref{sub.proj}, there exist $\xi^{k+1} \in \mathbb{R}_{+}$ and $\nu_{i}^{k+1} \in \mathbb{R}_{+}$ such that for $i \in \mathcal{I}^k$
\begin{subequations}\label{eq: Proj_kkt}
	\begin{alignat}{1}
		&\frac{1}{\beta}(x_{i}^{k+1} \!-\! x_{i}^{k}) \!+\! \nabla_i f(\bm{x}^{k}) + \xi^{k+1}\text{sgn}(x_i^{k+1}) w_{i}^{k} - \nu_{i}^{k+1} \textrm{sgn}(x_{i}^{k})\!=\! 0, i\in \I^{k}, \label{eq: Proj_kkt1}\\
		&\sum_{i\in \I^{k}} w_{i}^{k}(\vert x_{i}^{k+1}\vert - \vert x_{i}^{k}\vert)  \leq 0,\quad \bm{x}_{\mathcal{A}^k}^{k+1} = \bm{0}, \ \textrm{sgn}(\bm{x}^{k}_{\mathcal{I}^k}) \circ \bm{x}^{k+1}_{\mathcal{I}^k} \geq \bm{0}, \label{eq: Proj_kkt2}\\
		&\xi^{k+1}\Big( \sum_{i\in \I^{k}} w_{i}^{k}(\vert x_{i}^{k+1}\vert - \vert x_{i}^{k}\vert)\Big) =0, \ \nu^{k+1}_{i}\textrm{sgn}(x_{i}^{k})x_{i}^{k+1}=0,   \ i\in \I^k.  \label{eq: Proj_kkt3}
	\end{alignat}
\end{subequations}
We also use the subproblems solutions to characterize the optimal solutions of \eqref{eq: main_Opt}.
\begin{lemma}\label{Lem: subproblem_property}
	Consider \eqref{eq: main_Opt}. Given any $\bar{\bm{x}} \in \mathcal{B}_{\ell_{p}}$ and suppose $\nabla f(\bar{\bm{x}}) \neq \bm{0}$. The following hold:
	\begin{itemize}
		\item[(i)] If $\bar{\bm{s}} \in \argmin\limits_{\bm{s}\in\mathcal{B}_{\ell_p}} Q(\bm{s};\bar{\bm{x}})$ such that $\la \nabla f(\bar{\bm{x}}), \bar{\bm{s}} - \bar{\bm{x}}\ra = 0$, then $\bar{\bm{x}}$ is first-order stationary to \eqref{eq: main_Opt}.
		\item[(ii)] Let $\hat{\bm{x}} = \argmin\limits_{\bm{x}\in\mathcal{B}_{\ell_1^{\bar{\bm{w}}}}} P(\bm{x};\bar{\bm{x}})$. If $\hat{\bm{x}} = \bar{\bm{x}}$, then
		$\bar{\bm{x}}$ is first-order stationary to \eqref{eq: main_Opt}. Otherwise, $\hat{\bm{x}}  \in \textbf{int}\ \mathcal{B}_{\ell_p}$.
	\end{itemize}
\end{lemma}
\begin{proof}
	For (i),  it holds that $\bar{\bm{x}}$ is also an optimal solution of \eqref{eq: FW_linearization}, then $\bar{\bm{x}}$ satisfies \eqref{eq:optimality_FW_sub}. This implies that $\bar{\bm{x}}$ satisfies \eqref{Def SC} (by Theorem \ref{Theo: FW_reformulation}) and thus $\bar{\bm{x}}$ is first-order stationary to \eqref{eq: main_Opt}.
	
	For (ii), if $\bar{\bm{x}}$ is optimal to \eqref{sub.proj} and $\nabla f(\bar{\bm{x}}) \neq \bm{0}$, it then follows from $\Vert \bar{\bm{x}}\Vert_p^{p} = \gamma$ that \eqref{eq: Proj_kkt} reverts to \eqref{Def SC}, implying $\bar{\bm{x}}$ is first-order stationary to \eqref{eq: main_Opt}.  Otherwise, we know from \eqref{sub.proj} that $\hat{\bm{x}}$ and $\bar{\bm{x}}$ are in the same orthant.  Since 
	$\hat{\bm{x}} \ne \bar{\bm{x}} $ and recalling $\textrm{sgn}(\bar{\bm{x}}_{\mathcal{I}(\bar{\bm{x}})}) \circ \hat{\bm{x}}_{\mathcal{I}(\bar{\bm{x}})} \geq \bm{0}$ and $ \hat{\bm{x}}_{\mathcal{A}(\bar{\bm{x}})} = \bm{0}$, we have from \eqref{eq: well-posedness} and the strict concavity of $\Vert \cdot \Vert_{p}$ that
	$$ \Vert \hat\bx\Vert_p^p < \Vert \bar\bx\Vert_p^p + \sum_{i\in \I(\bar \bx)} p |\bar x_i|^{p-1} (\vert \hat x_{i}\vert - \vert \bar x_{i}\vert) \leq \gamma,$$
	completing the proof. 
\end{proof}

We present some useful results in the following lemma.
\begin{lemma}\label{finite.terminate} 
	Assume Algorithm \ref{alg.propose} does not terminate in finite iterations.  Then $\{k+1: k\in\mathcal{T}_2\} \subseteq  \mathcal{S}_2$ and  $\{k+1: k\in\mathcal{S}_2\} \subseteq \mathcal{S}_1$.  
	Therefore, one of the following cases holds. 
	\begin{itemize}
		\item[(i)] $|\mathcal{T}_1| = +\infty$,  $|\mathcal{T}_2| < +\infty$, $|\mathcal{S}_2| < +\infty$. 
		\item[(ii)]  $|\mathcal{S}_1 | = |\mathcal{S}_2| = +\infty$.
	\end{itemize}
\end{lemma}
\begin{proof} 
	If $k \in \mathcal{S}_2$ and the algorithm does not terminate at the $k+1$ iteration, then Lemma \ref{Lem: subproblem_property}(ii) implies that $\bx^{k+1} \in \textbf{int}\ \mathcal{B}_{\ell_{p}}$. Therefore,  subproblem \eqref{eq: FW_linearization} is solved at the $(k+1)$ iteration and hence $(k+1) \in \mathcal{S}_1$.  It follows that $|\mathcal{S}_2| = +\infty$ implies $|\mathcal{S}_1| = +\infty$.

	If $k\in \mathcal{T}_2$ and the algorithm does not terminate at the $k+1$ iteration, the definition of $\mathcal{T}_2$ implies that 
	$\bx^{k+1} \in \textbf{bd}\ \mathcal{B}_{\ell_{p}}$. Therefore, subproblem  \eqref{sub.proj} is solved at the $k+1$ iteration and 
	$(k+1) \in \mathcal{S}_2$.  It follows that $\vert \mathcal{T}_2 \vert = +\infty$ implies $\vert\mathcal{S}_2\vert = +\infty$.

	Overall, we have two cases based on $|\mathcal{S}_2|$.  (i) $|\mathcal{S}_2| = +\infty$ implies $|\mathcal{S}_1| = +\infty$. (ii) $|\mathcal{S}_2| < +\infty$. In the latter case, if assuming  $|\mathcal{T}_2| = +\infty$, then it implies $|\mathcal{S}_2| = +\infty$ which however causes a contradiction, thus   $|\mathcal{T}_2|<+\infty$ holds true, and it further  implies from $|\mathcal S_1|=+\infty$ that   $|\mathcal{T}_1|=+\infty$. 
\end{proof}

\subsection{Global Convergence}

We now analyze the convergence result 
when an infinite sequence $\{x^k\}$ is generated. We first summarize a useful property of sequence $\{\bx^{k+1}\}$ in the following lemma.
\begin{lemma}\label{Lemma: Basic_property_Proj}
	Consider Algorithm \ref{alg.propose} for solving \eqref{eq: main_Opt} under Assumption \ref{BasicAssum}. Let $\{\bm{x}^{k}\}$ be an infinite sequence generated by Algorithm \ref{alg.propose}. The following statements hold.
	\begin{itemize}
		\item[(i)] For $k \in \mathcal{S}_1 $, let $T_k$ be the total number of evaluations of the sufficient decrease condition implemented in Algorithm \ref{alg.subprocedure_fw} up to iteration $k\geq 1$. Then it holds for $k \in \mathcal{T}_1$ that
		$$T_{k} \leq \left[1 + \frac{\log\zeta}{\log\tau}\right](\vert \mathcal{T}_1 \vert+1) + \frac{1}{\log\tau}\max\left\lbrace\log\frac{\tau L_{f}}{L_{-1}},0 \right\rbrace.$$
		\item[(ii)] For $k \in \mathcal{S}_2 $, it holds that $\I(\bm{x}^{k}) \neq \emptyset$ and $\I(\bm{x}^{k+1}) \subseteq \I(\bm{x}^{k})$. 
	\end{itemize}
\end{lemma}
\begin{proof}
	The proof of statement (i) follows a similar way to that of \cite[Theorem 1]{pedregosa2020linearly}. For the sake of completeness, we here present the detailed proof. For each $i \in \mathcal{T}_1$, $t_i$ denotes the number of evaluations of the sufficient decrease condition at $i$th iteration. By Lines \ref{line:localLip} and \ref{line:tauMultiplication} in Algorithm \ref{alg.subprocedure_fw}, we know that the algorithm increases the estimation of $L_f$ by a power factor of $\tau>1$ whenever the sufficient decrease condition is not satisfied. Hence, we know that
	\begin{equation}\label{eq:LipNum}
		L_{i} = \zeta L_{i-1}\tau^{t_{i}-1}.
	\end{equation}
	Taking logarithms on both sides of \eqref{eq:LipNum} gives
	\begin{equation}\label{eq:LipNum1}
		t_i = 1 + \frac{\log\zeta}{\log\tau} + \frac{1}{\log\tau}\log\left(\frac{L_{i}}{L_{i-1}}\right).
	\end{equation}
	Summing up \eqref{eq:LipNum1} over $i \in  \mathcal{T}_1 $ up to  $k$th iteration yields
	\begin{equation}
		\begin{aligned}
			T_{k} \leq \sum_{i \in  \mathcal{T}_1, i\leq k }t_{i} &= \left[1 + \frac{\log\zeta}{\log\tau}\right](\vert \mathcal{T}_1 \vert+1) + \frac{1}{\log\tau}\log\left( \frac{L_{k}}{L_{-1}}\right)\\
			&\leq \left[1 + \frac{\log\zeta}{\log\tau}\right](\vert \mathcal{T}_1 \vert+1) + \frac{1}{\log\tau}\max\left\{\log\frac{\tau L_{f}}{L_{-1}},0\right\},
		\end{aligned}
	\end{equation}
	where the second inequality follows from $L_{k} \leq \max\{\tau L_{f}, L_{-1}\}, \forall k \in \mathbb{N}$ established in \cite[Proposition 2]{pedregosa2020linearly}, as desired.
	
	For (ii), it is trivial to see that $\I(\bm{x}^{k}) \neq \emptyset$ since $\Vert\bm{x}^{k}\Vert_{p}^{p} = \sum_{i \in \I(\bm{x}^{k})}\vert x_{i}^k\vert^{p} = \gamma$. Consequently, $w_{i}^{k}$ is finite for $i \in \I(\bm{x}^{k})$ and hence the interior of the constraint set of \eqref{sub.proj} is nonempty. Then, by the optimality of $\bm{x}^{k+1}$, we have that $x_i^{k+1} = 0$ if $x_{i}^{k} = 0$ for any $i \in \mathcal{A}^k$ by recalling the constraint $\bm{x}_{\mathcal{A}^k} = \bm{0}$ imposed in \eqref{sub.proj}, meaning $\mathcal{I}^{k+1} \subseteq \mathcal{I}^k$ for $k \in \mathcal{S}_2 \subseteq  \mathbb{N}$. 
 \end{proof}


The global convergence results are established in this subsection. We next show the sufficient reduction in $f$, $P$, and $Q$ in the following lemma, and its proof can be found in Appendix \ref{app_C}.
\begin{lemma}\label{Lem: Obj_decrease}
	Suppose $\{\bm{s}^{k}\}$ and $\{\bm{x}^{k}\}$ are generated by Algorithm \ref{alg.propose} with $\bm{x}^{0} \in \mathcal{B}_{\ell_p}$. Let $\tau > 1$. Then the following statements hold.
	\begin{itemize}
		\item[(i)] For $k\in \mathcal{T}_1$, $0 \leq \min\Big(\frac{\Delta Q^2(\bm{s}^{k};\bm{x}^{k})}{8\max\left(\tau L_{f}, L_{-1}\right) \gamma^{2/p}},\frac{\Delta Q(\bm{s}^{k};\bm{x}^{k})}{2} \Big)\leq \Delta f(\bm{x}^{k+1})$.\\ For $k\in\mathcal{T}_{2}$, $0 \leq
		\frac{\Delta Q^{2}(\bm{s}^{k};\bm{x}^{k})}{8\max\left(\tau L_{f}, L_{-1}\right) \gamma^{2/p}} \leq \Delta f(\bm{x}^{k+1})$.
		
		\item[(ii)] For $k \in \mathcal{S}_2$, $\Vert\bm{x}^{k} - \bm{x}^{k+1}\Vert_{2}^{2} \leq  \Delta P(\bm{x}^{k+1};\bm{x}^{k}) \leq \Delta f(\bm{x}^{k+1}).$
		
		\item[(iii)]  $\lim\limits_{\substack{k\in\mathcal{S}_1\\k\to+\infty}}\Delta Q(\bm{s}^{k};\bm{x}^{k}) = 0$ and $\lim\limits_{\substack{k\in\mathcal{S}_2\\k\to+\infty}}\Delta P(\bm{x}^{k+1};\bm{x}^{k}) = 0$. Hence, $\lim\limits_{\substack{k\in\mathcal{S}_2\\k\to+\infty}} \Vert\bm{x}^{k+1} - \bm{x}^{k}\Vert_2^2 = 0$.
		
		\item[(iv)] If $\lim\limits_{\substack{k\in\mathcal{T}_2\\k\to+\infty}}\alpha_{\textrm{bis}}^{k} = 0$, then $\lim\limits_{\substack{k\in\mathcal{T}_2\\k\to+\infty}} \Vert\bm{x}^{k+1} - \bm{x}^{k}\Vert_2^2 = 0$.
	\end{itemize}
\end{lemma}

Now we show the main convergence result in the following.
\begin{theorem}[Subsequential convergence]\label{eq: main_theo}
	Suppose $\{\bm{x}^{k}\}$ is an infinite sequence generated by Algorithm \ref{alg.propose} with $\bm{x}^{0} \in \mathcal{B}_{\ell_p}$. Then the sequence $\{\bm{x}^{k}\}$ is bounded. Moreover, every cluster point $\bm{x}^{*}$ of $\{\bm{x}^{k}\}$ is first-order stationary for \eqref{eq: main_Opt}.
\end{theorem}

\begin{proof}
	Since $\{\bm{x}^{k}\} \subset \mathcal{B}_{\ell_p}$ and $\mathcal{B}_{\ell_p}$ is bounded, $\{\bm{x}^{k}\}$ hence is bounded and has cluster points. To establish the global convergence results, we consider two cases. 
	
	Case (i):  $\bm{x}^{*}$ is  a cluster point of $\mathcal{S}_1$.   By Lemma \ref{Lem: subproblem_property}(i), it suffices to prove that $\bm{x}^{*}\in\argmin\limits_{\bm{s} \in \mathcal{B}_{\ell_p}}Q(\bm{s};\bm{x}^{*})$. We prove this by contradiction. Suppose this is not true. There exists $\tilde{\bm{x}}$ such that $\tilde{\bm{x}} \in \mathcal{B}_{\ell_p}$ and that $\epsilon := Q(\bm{x}^{*};\bm{x}^{*}) - Q(\tilde{\bm{x}};\bm{x}^{*}) > 0$. Consider a subsequence $\mathcal{U}$ with $\{\bm{x}^{k}\}_\mathcal{U} \to \bm{x}^{*}$. By Lemma \ref{Lem: Obj_decrease}(iii), for sufficiently large $k \in \mathcal{U} \subseteq \mathcal{T}_1$, it holds
	\begin{equation}\label{eq: FW_obj_decrease}
		0< Q(\bm{x}^{k};\bm{x}^{k}) - Q(\bm{s}^{k};\bm{x}^{k}) <\epsilon/2.
	\end{equation}
	By the continuity of $Q(\bm{s};\bm{x}^{k})$ with respect to $\bm{s}$ over $\mathcal{B}_{\ell_p}$, we can pick sufficiently large $k \in \mathcal{U}$ such that 
	\begin{equation}\label{eq: formula_2}
		\vert Q(\bm{x}^{k};\bm{x}^{k}) - Q(\bm{x}^{*};\bm{x}^{*})\vert <\epsilon/4\quad\text{ and }\quad \vert Q(\tilde{\bm{x}};\bm{x}^{k}) - Q(\tilde{\bm{x}};\bm{x}^{*})\vert < \epsilon/4.
	\end{equation}
	We have
	\begin{equation}\label{eq: dirive_contradiction}
		\begin{aligned}
			&\quad\ 
			Q(\bm{x}^{k};\bm{x}^{k}) - Q(\tilde{\bm{x}};\bm{x}^{k})\\
			&= Q(\bm{x}^{k};\bm{x}^{k}) - Q(\bm{x}^{*};\bm{x}^{*}) + Q(\bm{x}^{*};\bm{x}^{*}) -Q(\tilde{\bm{x}};\bm{x}^{*})+ Q(\tilde{\bm{x}};\bm{x}^{*})- Q(\tilde{\bm{x}};\bm{x}^{k}) \\
			&\geq-\vert Q(\bm{x}^{k};\bm{x}^{k}) - Q(\bm{x}^{*};\bm{x}^{*}) \vert + \vert Q(\bm{x}^{*};\bm{x}^{*}) -Q(\tilde{\bm{x}};\bm{x}^{*})\vert \vert Q(\tilde{\bm{x}};\bm{x}^{*})-Q(\tilde{\bm{x}};\bm{x}^{k})\vert\\
			&>-\epsilon/4 + \epsilon - \epsilon/4= \epsilon/2.
		\end{aligned}
	\end{equation}
	This contradicts the optimality of $\bm{s}^{k}$ for \eqref{eq: FW_linearization}. Hence, by Lemma \ref{Lem: subproblem_property}(i), $\bm{x}^{*}$ is first-order stationary to \eqref{eq: main_Opt}. This completes the proof. 
	
	Case (ii):  $\bm{x}^{*}$ is  a cluster point of $\mathcal{S}_2$.    By Lemma \ref{Lem: subproblem_property}(ii), it suffices to show that $\bm{x}^{*} = \argmin_{\bm{x}\in \mathcal{B}_{\ell_{1}}^*} P(\bm{x};\bm{x}^*)$. We prove this by contradiction. Suppose this is not true, and thus, there exists an optimal solution $\tilde{\bm{x}}$ such that $\sum_{i\in \I(\bm{x}^{*})} w_{i}^{k}(\vert\tilde{x}_{i}\vert - \vert x_{i}^{*}\vert) \leq 0$, $\tilde{\bm{x}}_{\mathcal{A}^{k}} = \bm{0}$, $\tilde{\bm{x}}_{\mathcal{I}^*} \circ \textrm{sgn}( \bm{x}^{*}_{\mathcal{I}^*}) \geq \bm{0}$ and that $\epsilon = P(\bm{x}^{*};\bm{x}^{*}) - P(\tilde{\bm{x}};\bm{x}^{*}) >0$. Consider a subsequence $\mathcal{U}$ with $\{\bm{x}^{k}\}_\mathcal{U} \to \bm{x}^{*}$. By Lemma \ref{Lem: Obj_decrease}(iii), there exists $k_1$ such that for any $k > k_1$ and $k \in \mathcal{U} \subseteq \mathcal{S}_2$
	\begin{equation}\label{eq: reduction}
		0 \leq P(\bm{x}^{k}; \bm{x}^{k}) - P(\bm{x}^{k+1}; \bm{x}^{k}) <\epsilon/4.
	\end{equation}
	Now consider the projection of $\tilde{\bm{x}}$ onto $\mathcal{B}_{\ell_{1}}^{k}$, denoted as $\tilde{\bm{x}}^{k}$. Since $P(\bm{x};\bm{x}^{k})$ is continuous with respect to $\bm{x}$ over $\mathcal{B}_{\ell_{1}}^{k}$ and $\{\bm{x}^{k}\}_\mathcal{U} \to \bm{x}^{*}$, there exists $k_2$ and $k_3$ such that for any $k > k_2$ and $k\in\mathcal{U}$,
	\begin{equation*}\label{eq: two_assum1}
		\vert P(\bm{x}^{k};\bm{x}^{k}) - P(\bm{x}^{*};\bm{x}^{*})\vert< \epsilon/4,
	\end{equation*}
	and that for any $k > k_3$ and $k \in \mathcal{U}$,
	\begin{equation*}\label{eq: two_assum2}
		\vert P(\tilde{\bm{x}}^{k};\bm{x}^{*}) - P(\tilde{\bm{x}};\bm{x}^{*})\vert< \epsilon/4.
	\end{equation*}
	Moreover, by the continuity of $P(\tilde{\bm{x}}^k;\bm{x})$ with respect to $\bm{x}$ over $\mathcal{B}_{\ell_{1}}^k$, we have
	\[
	\vert P(\tilde{\bm{x}}^k;\bm{x}^{*}) - 
	P(\tilde{\bm{x}}^{k};\bm{x}^{k})\vert < \epsilon/4.
	\]
	For $k > \max\{k_1, k_2, k_3\}$, it then follows that
	\begin{equation}\label{eq: contradiction}
		\begin{aligned}
			&\quad\ P(\bm{x}^{k};\bm{x}^{k}) - P(\tilde{\bm{x}}^{k};\bm{x}^{k})\\
			&= P(\bm{x}^{k};\bm{x}^{k}) - P(\bm{x}^{*};\bm{x}^{*}) + P(\bm{x}^{*};\bm{x}^{*})
			- P(\tilde{\bm{x}};\bm{x}^{*})+ P(\tilde{\bm{x}};\bm{x}^{*}) -P(\tilde{\bm{x}}^{k};\bm{x}^{*})\\
			&\quad
			+P(\tilde{\bm{x}}^{k};\bm{x}^{*}) - 
			P(\tilde{\bm{x}}^{k};\bm{x}^{k})\\
			&\geq -\vert P(\bm{x}^{k};\bm{x}^{k}) - P(\bm{x}^{*};\bm{x}^{*})\vert + \vert P(\bm{x}^{*};\bm{x}^{*})
			- P(\tilde{\bm{x}};\bm{x}^{*})\vert - \vert P(\tilde{\bm{x}};\bm{x}^{*}) - P(\tilde{\bm{x}}^{k};\bm{x}^{k})\vert\\
			&\quad - \vert P(\tilde{\bm{x}}^{k};\bm{x}^{*}) - 
			P(\tilde{\bm{x}}^{k};\bm{x}^{k})\vert\\
			&>-\epsilon/4 + \epsilon -\epsilon/4 -\epsilon/4 = \epsilon/4,
		\end{aligned}
	\end{equation}
	contradicting \eqref{eq: reduction}. This indicates that $\tilde{\bm{x}}^{k}$ is a feasible iterate for the $k$th subproblem \eqref{sub.proj} but with lower objective value than $\bm{x}^{k+1}$ and hence contradicts the optimality of $\bm{x}^{k+1}$ to \eqref{eq: FW_linearization}. Therefore, $\bm{x}^{*}$ is first-order stationary to \eqref{eq: main_Opt} by Lemma \ref{Lem: subproblem_property}(ii). 
\end{proof}
\subsection{Convergence Rate Analysis}\label{Sec: complexity}
In this subsection, we analyze the local convergence rate of Algorithm \ref{alg.propose}. We first define the following convergence criteria to measure the optimality errors at $\bm{x}^{k}$, which  correspondingly associate with the subsequences $\mathcal{S}_1$ and $\mathcal{S}_2$:
\begin{subequations}\label{eq: residuals1}
	\begin{alignat}{1}
		&E_{\mathcal{S}_1}(\bm{x}^{k}) := Q(\bm{x}^{k};\bm{x}^{k}) - Q(\bm{s}^{k};\bm{x}^{k}), \label{eq: residua2}\\
		&E_{\mathcal{S}_2}(\bm{x}^{k},\xi^{k+1}) := \Vert \nabla_{\mathcal{I}^k} f(\bm{x}^{k}) \circ \bm{x}^{k+1}_{\mathcal{I}^k} + \xi^{k+1} p \vert\bm{x}^{k}_{\mathcal{I}^k}\vert^{p-1} \circ \text{sgn}(\bm{x}_{\mathcal{I}^k}^k) \circ \bm{x}^{k+1}_{\mathcal{I}^k}\Vert_{1}^2. \label{eq: residua3}
	\end{alignat}
\end{subequations}
For $\mathcal{S}_1$, we follow \cite{pedregosa2020linearly} to consider $E_{\mathcal{S}_1}(\bm{x}^{k})$ as the so-called \textit{duality gap} $\la \nabla f(\bm{x}^{k}), \bm{x}^{k} - \bm{s}^{k}\ra$, a standard practice within the framework of the Frank-Wolfe method. Note that it is equivalent to the gap in the objective value of the FW subproblem. From  Lemma \ref{Lem: subproblem_property}(i),  $E_{\mathcal{S}_1}(\bm{x})$ vanishes only at a stationary point. As for $\mathcal{S}_2$, the optimality residual is indicated by the scaled optimality conditions of Corollary \ref{LocalImpliesSP}.

The following theorem establishes an upper bound of the iteration number when one of the defined optimality errors is less than the specified error tolerance. Specifically, it indicates that at most $O(1/\epsilon^2)$ iterations, one can have $E_{\mathcal{S}_1}^{2}(\bm{x}^{k})\leq\epsilon$ on the subsequence $\mathcal{S}_1$, or $E_{\mathcal{S}_2}(\bm{x}^{k})\leq\epsilon$ on the subsequence $\mathcal{S}_2$.
\begin{theorem}\label{Lem: complexity1}
	Let $\epsilon > 0$ be given. Suppose the sequence $\{\bm{x}^{k}\}$ is generated by Algorithm \ref{alg.propose} with $\bm{x}^{0} \in \mathcal{B}_{\ell_{p}}$ and clusters at $\bm{x}^{*} \in \mathcal{B}_{\ell_p}$. The the maximal number of iterations of Algorithm \ref{alg.propose} to attain $$\min\{E_{\mathcal{S}_1}^{2}(\bm{x}^{k}), E_{\mathcal{S}_2}(\bm{x}^{k},\xi^{k+1})\} \leq \epsilon$$
	is $\frac{f(\bm{x}^0) - \underline{f}}{\min\{C_1, C_2\}\epsilon}$, where $\underline{f} = \lim\limits_{j \to +\infty} f(\bm{x}^{k_{j}}) = f(\bm{x}^{*})$, $C_1:=\tfrac{1}{8\max\left(\tau L_{f}, L_{-1}\right) \gamma^{2/p}}$ and $C_2 = \frac{\beta^{2}}{n \gamma^{2/p}}$.
\end{theorem}
\begin{proof}
	For $k \in \mathcal{S}_1$,  since we aim to evaluate the worst-case iteration complexity of Algorithm \ref{alg.propose}, by Lemma \ref{Lem: Obj_decrease}(i), we have 
	\begin{equation}\label{eq: complexity_ine2_2}
		\begin{aligned}
			E^{2}_{\mathcal{S}_1}(\bm{x}^{k})
			=   \Delta Q^{2}(\bm{s}^{k};\bm{x}^{k})
			\leq 8\max\left(\tau L_{f}, L_{-1}\right) \gamma^{2/p}\Delta f(\bm{x}^{k+1}). 
		\end{aligned}
	\end{equation}
	Therefore, we have
	\begin{equation}\label{eq: Residual_FW}
		\frac{1}{8\max\left(\tau L_{f}, L_{-1}\right) \gamma^{2/p}} E^{2}_{\mathcal{S}_1}(\bm{x}^{k}) \leq  \Delta f(\bm{x}^{k+1}).
	\end{equation}
	
	By \eqref{eq: Proj_kkt1}, we have for $k \in \mathcal{S}_2$ and $i \in \I^k$
	\begin{equation}\label{eq: Ine_useful}
		\nabla_{i} f(\bm{x}^{k}) + \xi^{k+1}w_{i}^{k}\text{sgn}(x_i^k)
		= \frac{1}{\beta}(x_{i}^{k} - x_{i}^{k+1}) + \nu_i^{k+1} \text{sgn}(x_{i}^{k}).
	\end{equation}
	Multiplying $x_i^{k+1}$ for $i \in \mathcal{I}^k$ on both sides of \eqref{eq: Ine_useful}, we obtain
	\begin{equation}\label{eq: Ine_useful1}
		\begin{aligned}
			&\quad\  \nabla_{i} f(\bm{x}^{k}) x_i^{k+1} + p\xi^{k+1}\text{sgn}(x_i^k)\vert x_{i}^{k} \vert^{p-1} x_i^{k+1}\\
			&= \frac{1}{\beta}(x_{i}^{k} - x_{i}^{k+1}) x_i^{k+1} + \nu_i^{k+1} \text{sgn}(x_{i}^{k})x_i^{k+1} = \frac{1}{\beta}(x_{i}^{k} - x_{i}^{k+1}) x_i^{k+1},
		\end{aligned}
	\end{equation}
	where the last equality is true because $\nu_i^{k+1} = 0$ holds due to the complementary condition \eqref{eq: Proj_kkt3}.
	Recall \eqref{eq: residua3}, we have
	\begin{equation}\label{eq: complexity_ine1}
		\begin{aligned}
			\ E_{\mathcal{S}_2}(\bm{x}^{k},\xi^{k+1}) &=\Vert p\xi^{k+1} \text{sgn}(\bm{x}^{k}_{\I^k})\circ\vert\bm{x}^{k}_{\I^k}\vert^{p-1} \circ \bm{x}_{\mathcal{I}^{k}}^{k+1} + \nabla_{\I^k} f(\bm{x}^{k})\circ \bm{x}_{\mathcal{I}^{k}}^{k+1} \Vert_{1}^2\\
			&= \frac{1}{\beta^2}\Vert(\bm{x}^{k}_{\I^k} - \bm{x}^{k+1}_{\I^k})\circ \bm{x}_{\mathcal{I}^{k}}^{k+1}\Vert_{1}^2 \leq \frac{1}{\beta^2}\Vert\bm{x}^{k}_{\I^k} - \bm{x}^{k+1}_{\I^k}\Vert_{1}^2 \Vert\bm{x}_{\mathcal{I}^{k}}^{k+1}\Vert_{1}^2\\
			& \leq \frac{\vert\mathcal{I}^{k}\vert}{\beta^2}\Vert\bm{x}^{k}_{\I^k} - \bm{x}^{k+1}_{\I^k}\Vert_{2}^2 \Vert\bm{x}_{\mathcal{I}^{k}}^{k+1}\Vert_{1}^2\leq \frac{\vert\mathcal{I}^{k}\vert \gamma^{\frac{2}{p}}}{\beta^2}\Vert\bm{x}^{k}_{\I^k} - \bm{x}^{k+1}_{\I^k}\Vert_{2}^2\\
			&\leq \frac{\vert\mathcal{I}^{k}\vert \gamma^{\frac{2}{p}}}{\beta^2}\Delta P(\bm{x}^{k+1};\bm{x}^{k}) \leq \frac{n \gamma^{\frac{2}{p}}}{\beta^2}\Delta f(\bm{x}^{k+1}),
		\end{aligned}
	\end{equation}
	where the second equality holds due to \eqref{eq: Ine_useful1}, the second inequality is right because of $\Vert \bm{x} \Vert_{1} \leq \Vert \bm{x} \Vert_2$, $\forall \bm{x} \in \mathbb{R}^{n}$, the third inequality holds because $\Vert\bm{x}^{k}\Vert_p \leq 
	\gamma^{1/p}$ implies $\Vert\bm{x}^{k}\Vert_{1} \leq \gamma^{1/p}$, and the fourth and the last inequality hold due to Lemma \ref{Lem: Obj_decrease}(ii).  
	
	Rearranged, \eqref{eq: complexity_ine1} gives
	\begin{equation}\label{eq: Residual_Proj}
		\frac{\beta^{2}}{n \gamma^{2/p}}E_{\mathcal{S}_2}(\bm{x}^{k},\xi^{k+1}) \leq \Delta f(\bm{x}^{k+1}).
	\end{equation}
	
	Consider now $\min\{E^{2}_{\mathcal{S}_1}(\bm{x}^{k}), E_{\mathcal{S}_2}(\bm{x}^{k}, \xi^{k+1})\} > \epsilon$.  Let $C_1:=\tfrac{1}{8\max\left(\tau L_{f}, L_{-1}\right) \gamma^{2/p}}$ and $C_2 = \frac{\beta^{2}}{n\gamma^{2/p} }$.  Summing \eqref{eq: Residual_FW} and \eqref{eq: Residual_Proj} both sides from $t = 0,\ldots,k$ yields
	\begin{equation*}
		\begin{aligned}
			\sum_{t\in\mathcal{S}_1\cap[k]} C_1E^{2}_{\mathcal{S}_1}(\bm{x}^{t}) +  \sum_{t\in\mathcal{S}_2\cap[k]}C_2 E_{\mathcal{S}_2}(\bm{x}^{t}, \xi^{t+1}) &\leq \sum_{t=0}^k \Delta f(\bm{x}^{t+1}) \\&= f(\bm{x}^0) - f(\bm{x}^{k+1})
			\leq f(\bm{x}^0) - \underline{f}.
		\end{aligned}
	\end{equation*}
	Rearranged, we have
	
	\begin{equation*}
		\begin{aligned}
			f(\bm{x}^0) - \underline{f} 
			&\ge \sum_{t\in\mathcal{S}_1\cap[k]} C_1E^{2}_{\mathcal{S}_1}(\bm{x}^{t}) + \sum_{t\in\mathcal{S}_2\cap[k]}C_2 E_{\mathcal{S}_2}(\bm{x}^{t}, \xi^{t+1}) \\
			&\geq \sum_{t\in\mathcal{S}_1\cap[k]} C_1E^{2}_{\mathcal{T}_1}(\bm{x}^{t})  + \sum_{t\in\mathcal{S}_2\cap[k]}C_2 E_{\mathcal{S}_2}(\bm{x}^{t}, \xi^{t+1})\\
			&\geq \min\{C_1, C_2\} (\vert\mathcal{S}_1\cap[k]\vert\epsilon  +\vert\mathcal{S}_2\cap[k]\vert\epsilon) \\
			&= \min\{C_1, C_2\} (\vert\mathcal{S}_1\cap[k]\vert + \vert\mathcal{S}_2\cap[k]\vert)\epsilon\\
			&= \min\{C_1, C_2\}k\epsilon,
		\end{aligned}
	\end{equation*}
	implying
	\begin{equation}
		k \leq \frac{f(\bm{x}^0) - \underline{f}}{\min\{C_1, C_2\}\epsilon}.
	\end{equation}
	This establishes an upper bound for the maximal iteration number $k$ for satisfying $\max\{E^{2}_{\mathcal{S}_1}(\bm{x}^{k}), E_{\mathcal{S}_2}(\bm{x}^{k}, \xi^{k+1})\} < \epsilon$.
	\end{proof}

\section{Numerical Experiments}\label{Exp_sec}
In this section, we conduct several numerical experiments both on synthetic data and real-world data to evaluate the performance of the proposed Algorithm \ref{alg.propose}. In the first experiment, we deliver numerical comparisons of our algorithm with an MM-type method \cite{yang2022towards} and two DC-type methods \cite{boob2020feasible,zeng2024frankwolfetype}, through testing the $\ell_p$-ball projection problem \cite{yang2022towards}. The second experiment is to test the proposed algorithm on sparse signal recovery problem \cite{oymak2017sharp}, which aims to recover sparse signals from linear measurements. In addition, we test the proposed algorithm on real-world image reconstruction problems. All codes are implemented in Python\footnote[1]{ \url{https://github.com/Optimizater/Hybrid-1st-order-Algorithm}} and run on a laptop with Intel Core CPU i9-13900K at 3.0GHz and 64GB of main memory. 

To terminate the root-finding method for solving $\chi(\alpha)=0$  in the FW step,  we deem the bisection method is successful  if it finds an $\alpha$ satisfying 
\begin{equation}\label{eq: root-err tol}
	0< \vert \chi(\alpha)\vert  < \delta_{\textrm{bis}}.
\end{equation}
Meanwhile,  we also use this criterion in our implementation to determine whether the current point is on the boundary.  Specifically,  we deem an iterate $\bm{x}^{k}$ lies on the boundary of the $\ell_{p}$ ball if $\vert\Vert\bm{x}^{k}\Vert_{p}^{p} -\gamma\vert \leq \delta_{\textrm{bis}} $ holds; otherwise, it is deemed in the interior of $\mathcal{B}_{\ell_p}$. In all experiments, $\delta_{\textrm{bis}}=10^{-10}$ is adopted. 
We claim that this practical change does not affect the convergence result since   
\begin{itemize}
	\item[(i)] criterion \eqref{eq: root-err tol} guarantees that the proposed algorithm is still well-defined in the sense that all iterates are feasible.
	\item[(ii)] criterion \eqref{eq: root-err tol}  is used to trigger the gradient projection subproblem.  As analyzed in \S\ref{Sec: complexity}, this user-specified constant $\delta_{\textrm{bis}}$ generally does not affect the complexity analysis results if we choose $\delta_{\textrm{bis}}\ll \epsilon$ since it does not contribute much in the optimality error of $\bm{x}^{k+1}$ and therefore can be ignored. 
\end{itemize}
In all tests, we employ the weighted generalization of the sort-based algorithm\cite[Algorithm 2]{perez2022efficient} as the subproblem solver for \eqref{sub.proj}. We terminate the algorithm if one of the following stopping criteria is met:
\begin{itemize}
	\item[(i)]  $\Vert\bm{x}^{k+1} - \bm{x}^k\Vert_2 < 10^{-8} $ if $\bm{x}^{k} \in \textbf{bd}\ \mathcal{B}_{\ell_p}$;
	
	\item[(ii)]  $ \la \nabla f(\bm{x}^k), \bm{s}^k - \bm{x}^k \ra  < 10^{-8} $ if $\bm{x}^{k} \in \textbf{int}\ \mathcal{B}_{\ell_p}$.
\end{itemize}
\subsection{Euclidean Projection onto the $\ell_p$ Ball}
If the function $f$ in problem \eqref{eq: main_Opt} represents a squared loss function, then the Euclidean projection onto the $\ell_p$ ball can be formulated as 
\begin{equation}\label{prob:proj}
	\begin{aligned}
		\min_{\bm{x}\in\mathbb{R}^{n}}&\quad f(\bm{x}) := \frac{1}{2} \Vert\bm{x} - \bm{y}\Vert_2^2  \\
		\text{s.t.} & \quad \Vert\bm{x}\Vert_p^p \leq \gamma,
	\end{aligned} 
\end{equation}
where $ \bm{y} \in \mathbb{R}^n $ is a given point to be projected and $ \gamma \in \mathbb{R}_{++} $ is the radius of the $\ell_{p}$ ball. In particular, we vacuously assume $\Vert\bm{y}\Vert_{p}^{p} > \gamma$. 

In this test, we deliver a set of performance comparisons of the proposed algorithm against several state-of-the-art methods. These include an iteratively reweighted $ \ell_1 $-ball projection algorithm (IRBP) introduced in \cite{yang2022towards}, a level-constrained proximal point (LCPP) method presented by \cite{boob2020feasible} and a Frank-Wolfe type algorithm outlined in the work of \cite{zeng2024frankwolfetype}. For simplicity, we refer to the Frank-Wolfe type algorithm as ``FW-Zeng". Notably, it has been documented that the effectiveness of LCPP and FW-Zeng in solving problem \eqref{eq: main_Opt_General}. This naturally extends to encompass \eqref{prob:proj} as a special case. Below, We provide the following concise descriptions of each algorithm to facilitate a better understanding for our readers.
\begin{itemize}
	\item[(i)] IRBP is formulated within the MM algorithmic framework \cite{lange2016mm}. This algorithm introduces a smoothing parameter in the $\ell_{p}$ norm to have a continuously differentiable constraint function. IRBP operates on such a perturbed constrained function directly to form a sequence of convex weighted $\ell_1$ ball projection subproblems. To guarantee the convergence of IRBP, a scheduling strategy for adaptively reducing this smoothing parameter is developed. 
	
	\item[(ii)] LCPP is introduced within the DC algorithmic framework. Similar to IRBP, LCPP first introduces a smoothing parameter in the $\ell_{p}$ norm, then it expresses such perturbed $\ell_{p}$ norm as the difference of two convex functions. Consequently, LCPP generates a sequence of convex projection-type subproblems that differ from those of IRBP. However, it is noteworthy that LCPP does not explicitly provide a reduction rule for the introduced smoothing parameters. 
	
	\item[(iii)] FW-Zeng is proposed within the DC algorithmic framework. FW-Zeng can employ the same DC form of the perturbed constraint function as LCPP. Diverging from IRBP and LCPP, this algorithm incorporates a linearized objective function, resulting in the formation of a linear optimization subproblem. Similar to LCPP, FW-Zeng requires an updating rule for the introduced smoothing parameters.
\end{itemize}

We now specify the experimental setup for the aforementioned four algorithms. We generate $ \bm{y}  \sim\mathcal{N}(\bm{0},\bm{1}) $ and set the radius as $\gamma = 10^{-2} \Vert\bm{y}\Vert_p^p$. The initialization is performed with $\bm{x}^{0} = 0.3\gamma^{1/p}\vert\bm{y}\vert/\Vert\bm{y}\Vert_{p}$. It is apparent that the Lipschitz constant of the gradient of the objective in \eqref{prob:proj} is $L_f = 1$. In addition, 
\begin{itemize}
	\item[(i)] for the proposed algorithm, we set the step-size $\beta = 0.3/L_f = 0.3$ in Line \ref{line:gradStep} of Algorithm \ref{alg.propose}.
	
	\item[(ii)] For IRBP, the smoothing parameter is initialized as $\bm{\epsilon}^{0} = 0.6\gamma^{1/p}\vert\bm{y}\vert/\Vert\bm{y}\Vert_{p}$ to guarantee $\Vert\bm{x}^{0} + \bm{\epsilon}^{0}\Vert_{p}^{p}\leq \gamma$. Within the IRBP framework, two residuals, one associated with optimality and the other with feasibility, are defined to evaluate the solution's quality:    
	\begin{equation*}
		\begin{aligned}
			\text{R}_{\text{opt}}(\bm{x},\xi) = \frac{1}{n}\sum_{i=1}^n \left\vert(x_i - y_i)x_i + \xi p \vert x_i\vert^p \right\vert  \quad\textrm{ and }\quad
			\text{R}_{\text{fea}}(\bm{x}) = \vert \Vert\bm{x}\Vert_p^p - \gamma \vert.
		\end{aligned}
	\end{equation*}
	We terminate IRBP if $ \max\{\text{R}_{\text{opt}}, \text{R}_{\text{fea}}\} \leq 10^{-7} $.
	
	\item[(iii)] For LCPP, the smoothing parameter is initialized as $\epsilon^{0} = 0.8((\gamma - \Vert\bm{x}^{0}\Vert_p^{p})/n)^{1/p}$ to ensure $\Vert\bm{x}^{0} + \epsilon^{0}\Vert_{p}^{p}\leq \gamma$ by following \cite{boob2020feasible}. Following \cite[Table 2]{boob2020feasible}, the DC form of $\ell_{p}$ norm is defined as $g(\bm{x}) = p\epsilon^{p-1} \Vert \bm{x} \Vert_1 - h(\bm{x})$, where $ h(\bm{x}) = p\epsilon^{p-1} \Vert \bm{x} \Vert_1 - \Vert \bm{x} + \epsilon \Vert_{p}^{p}$. Inspired by \cite{chartrand2008iteratively,liu2022improved} and supported by the effectiveness through our numerical observation, the update of the smoothing parameter of LCPP is adopted according to the following rule, namely,
	\begin{equation}
		\begin{aligned}
			\left\{\begin{array}{ll}
				\epsilon^{k+1} = \epsilon^k/10,  & \textrm{ if }\ \Vert \bm{x}^{k} - \bm{x}^{k-1} \Vert_2 \leq \sqrt{\epsilon^k}; \\
				\epsilon^{k+1} = \epsilon^{k}, &  \textrm{ otherwise. }
			\end{array}\right.
		\end{aligned}
	\end{equation}
	Parameters such as the proximity parameter and initial radius (notion in the original text is $\gamma$ and $\eta_0$, respectively. Please refer to \cite[Algorithm 1]{boob2020feasible}) are set as suggested by \cite{boob2020feasible}\footnote{Based on our numerical observations, it appears that the increment in the radius constraint (referred to as $\delta_{k}$ in the original text) within the subproblem impacts the algorithm's performance in solving \eqref{prob:proj}. As proposed  in \cite{boob2020feasible}, we adopt the formula $\delta_{k} = \tfrac{\eta - \eta_0}{k(k+1)}$. It's worth noting that this increment, denoted as $\delta_{k}$, diminishes as $k$ increases. This reduction in $\delta_{k}$ may potentially impede the subproblem's radius $\eta_k$ from converging to the desired radius $\eta$, thereby affecting the feasibility numerically. To enhance the performance, our implementation involves truncating $\epsilon$ to $0$ and setting the subproblem's radius to $\eta$ when $\eta_k$ and $\eta$ are in close proximity. Additionally, we empirically establish a truncation threshold of $0.1$ for $p \in \{0.5, 0.7, 0.9\}$, and $4$ for $p = 0.3$. This configuration has demonstrated enhanced performance in our numerical experiments.}. The termination criterion for its subproblem is configured to trigger if $\Vert \bm{x}^{k} - \bm{x}^{k+1} \Vert_2 / (\Vert \bm{x}^{k+1}\Vert_2 + 10^{-10}) < 10^{-10}$ or if it exceeds $10$ iterations, as suggested by \cite{boob2020feasible}. At the outer loop level, termination occurs when $\Vert \bm{x}^{k} - \bm{x}^{k+1} \Vert_2 / (\Vert \bm{x}^{k+1}\Vert_2 + 1) < 10^{-8}$ or when the maximum number of iterations reaches $10^5$.
	\item[(iii)] The implementation of the FW-Zeng algorithm is publicly available at this link\footnote[1]{\url{https://github.com/zengliaoyuan/nonconvex_FW_code}}. In our test, FW-Zeng employs the same perturbation strategy and DC form of the $\ell_{p}$ norm constraint as utilized in LCPP. In addition, we set the parameters as $c = 10^{-3}, \eta = \frac{1}{2}, \alpha_0^0 = 1$ (as outlined in \cite[Algorithm 2]{zeng2024frankwolfetype})\footnote[2]{In our numerical observations applying FW-Zeng to address \eqref{prob:proj}, we observe that the algorithm generally fails to produce a satisfactory output within a two-minute time frame. This suggests the need for additional iterations to enhance its feasibility, even as the objective function exhibits minimal changes. It is noteworthy that achieving close-to-zero feasibility is challenging, particularly when dealing with stringent values of $p$ (e.g., $p < 0.5$). To enhance the algorithm's efficacy in solving \eqref{prob:proj}, we implement a strategy to truncate the perturbation to $0$ when the time limit approaches $200$ seconds.}. The termination criterion for FW-Zeng is set to activate when either  $\vert \langle \nabla f(\bm{x}^k), \bm{d}^k\rangle\vert / \vert f(\bm{x}^k)\vert  < 10^{-5}$ or $\Vert \bm{x}^{k} - \bm{x}^{k+1} \Vert_2 / (\Vert \bm{x}^{k+1}\Vert_2 + 1) < 10^{-8}$.
\end{itemize}

Numerical performance comparisons among the four methods are outlined in Table \ref{table.proj3}. A direct observation from Table \ref{table.proj3} reveals the general superiority of the proposed algorithm over comparative methods across all performance metrics considered for various $p$ values. Notably, the proposed algorithm exhibits the highest computational efficiency in all $p$ scenarios, concurrently demonstrating commendable optimality and feasibility. 

\begin{table}[htbp]
	\centering
	\caption{The performance comparison of considered algorithms for different $p$ values with $n = 10^5$.}
	\begin{tabular}{clllll}
		\toprule
		Value of $p$  & Algorithm & Obj. Value & $\textrm{R}_{\textrm{opt}}$  & $\textrm{R}_{\textrm{fea}}$  & Time(s) \\
		\midrule
		\multirow{4}[2]{*}{0.1} & IRBP  & 48080.25  & 3.08E-04 & 536.30  & 159.50  \\
		& LCPP  & 49942.32  & 6.78E-35 & 943.85  & 0.14  \\
		& FW-Zeng & 49942.32  & 1.22E-25 & 943.58  & 0.01  \\
		& Proposed & 46252.78  & 1.57E-08 & 1.03E-03 & 2.94  \\
		\midrule
		\multirow{4}[2]{*}{0.3} & IRBP  & 47766.00  & 3.80E-03 & 222.94  & 200.19  \\
		& LCPP  & 47126.27  & 3.57E-10 & 3.13E-06 & 65.44  \\
		& FW-Zeng & 48667.32  & 5.08E-03 & 503.68  & 170.13  \\
		& Proposed & 47104.16  & 2.38E-12 & 1.55E-07 & 2.49  \\
		\midrule
		\multirow{4}[2]{*}{0.5} & IRBP  & 47691.94  & 5.11E-09 & 6.59E-06 & 21.98  \\
		& LCPP  & 47689.57  & 3.06E-09 & 3.70E-09 & 87.41  \\
		& FW-Zeng & 48520.52  & 1.03E-04 & 188.44  & 200.00  \\
		& Proposed & 47637.75  & 4.85E-14 & 4.70E-08 & 0.86  \\
		\midrule
		\multirow{4}[2]{*}{0.7} & IRBP  & 48055.05  & 2.98E-10 & 6.53E-06 & 4.41  \\
		& LCPP  & 48054.02  & 1.29E-10 & 9.65E-11 & 72.89  \\
		& FW-Zeng & 48323.61  & 3.75E-04 & 73.41  & 200.00  \\
		& Proposed & 47827.22  & 3.36E-11 & 9.14E-10 & 1.75  \\
		\midrule
		\multirow{4}[2]{*}{0.9} & IRBP  & 47995.14  & 2.53E-12 & 2.61E-13 & 3.75  \\
		& LCPP  & 47947.68  & 7.77E-11 & 1.07E-12 & 53.75  \\
		& FW-Zeng & 48150.89  & 8.59E-03 & 6.52  & 200.00  \\
		& Proposed & 47960.94  & 2.50E-13 & 2.31E-12 & 2.34  \\
		\bottomrule
	\end{tabular}%
	\label{table.proj3}%
\end{table}%

\subsection{Sparse Signal Recovery}\label{Sec_sparsesignal}
In this subsection, we deliver a set of numerical experiments to evaluate the performance of the proposed algorithm for solving the sparse signal recovery problem. Given the measurement matrix $ \bm{A} \in \mathbb{R}^{m \times n} $ (often $m \leq n$ is assumed) and the observation vector $\bm{b} \in \mathbb{R}^m $. We consider the following $\ell_{p}$ ball constrained optimization problem
\begin{equation}\label{prob:sparse}
	\begin{aligned}
		\min_{\bm{x}\in\mathbb{R}^{n}}&\quad f(x) = \sum_{i=1}^{m}\mathcal{L}_i(\bm{x})  \\
		\text{s.t.} & \quad \Vert\bm{x}\Vert_p^p \leq \gamma.
	\end{aligned} 
\end{equation}
where $\mathcal{L}_i(\bm{x})$ is the loss function. In particular, we consider both the squared loss $\mathcal{L}_i(\bm{x}) = \tfrac{1}{2} \sum_{i=1}^{m}(\bm{a}_{i}^{T}\bm{x}-b_i)^2$ and the nonconvex Cauchy loss $\mathcal{L}_i(\bm{x}) = \sum_{i=1}^{m}\log(\tfrac{1}{2}(\bm{a}_{i}^{T}\bm{x}-b_i)^{2} + 1)$ \cite{carrillo2016robust}, where $\bm{a}_{i} \in \mathbb{R}^{n}$ is the $i$th row of $\bm{A}$.

The experiment is performed using synthetic data, where $ \hat{\bm{x}} $ represents the original $s$-sparse signal with $ s $ nonzero entries to be estimated. In all tests, we fix $ n = 1000 $ and set $ s = 100 $.  The number of measurements $m$ ranges from $50$ to $1000$. For each specified $m$, we repeat the following procedures $50$ times: 
\begin{itemize}
	\item[(i)]  Generate $ \hat{\bm{x}} $ by randomly selecting $ n-s $ entries and setting them to zero. Each nonzero entry is assigned a value of $-1$ or $+1$ with equal probability.
	\item[(ii)]  Generate the matrix $\bm{A}$ with entries drawn from a standard normal distribution and form $ \bm{b} = \bm{A}\hat{\bm{x}} + \bm{\epsilon} $ with $ \bm{\epsilon}_i \sim \mathcal{N}(0, 10^{-4}) $, $i \in [n]$.
	\item[(iii)] Resolve the associated optimization problems to obtain the estimated signal $ \bm{x}^{*} $.
\end{itemize}
As suggested by \cite{oymak2017sharp}, we declare a successful recovery for $\hat{\bm{x}}$ when  $ \Vert\bm{x}^{*} - \hat{\bm{x}}\Vert_2/\Vert\hat{\bm{x}}\Vert_2 < 10^{-3}$. Throughout all tests, we set $\gamma = s$.

\textbf{The convex squared loss.}\ In this experiment, we evaluate different algorithms with the squared loss. To facilitate performance comparison, the following algorithms are employed as benchmarks:
\begin{itemize}
	\item[(i)]  $\ell_{p}$-projected gradient descent ($\ell_{p}$-PGD) method is employed to solve \eqref{prob:sparse}. This algorithm iteratively solves a $\ell_{p}$ ball projection subproblem for $p \in [0,1]$, yielding an approximate solution \cite{bahmani2013unifying}. We consider $p$ values within the set $\{0, 0.5, 1\}$ in our tests.
	For $p = 0$, the algorithm aligns with the iterative hard thresholding (IHT) algorithm \cite{blumensath2009iterative}, designed for solving the $\ell_0$ ball-constrained problem. 
	For $p = 1$, the algorithm coincides with the $\ell_{1}$-PGD algorithm as presented in \cite{zhang2020inexact}. 
	For $0<p<1$, it is important to note that this algorithm is implemented without a theoretical guarantee for convergence \cite{perez2022efficient}.  
	\item[(ii)]  Proximal gradient method \cite{xu2012l_} for solving the $\ell_{p}$-regularized sparse signal recovery problem. Such an optimization problem is of the following form
	\begin{equation}\label{prob:sparse2}
		\min_{\bm{x}\in\mathbb{R}^{n}} \frac{1}{2} \Vert\bm{A}\bm{x} - \bm{b}\Vert_{2}^{2} + \lambda\Vert\bm{x}\Vert_{p}^{p},
	\end{equation}
	where $\lambda > 0$ is the so-called regularization parameter.  
\end{itemize} 

We initialize $\bm{x}^{0} = 0.9(\frac{\gamma}{\Vert\nu\Vert_{1}}\bm{\nu})^{1/p}$, where each entry of $\bm{\nu}$ is uniformly sampled from $\left[0,1 \right]$. In the proximal gradient method, $\tau$ is set to $10^{-4}$ through careful parameter tuning. It's worth noting that the Lipschitz constant $ L_f $ of the objective is equivalent to the largest eigenvalue $\lambda_{\max}$ of the positive semidefinite matrix $ \bm{A^T}\bm{A} $. We estimate it approximately using the power method \cite{journee2010generalized}. For all algorithms, the stepsize is uniformly set to $ \beta = 1/ \lambda_{\max} $ in the gradient projection step.

\begin{figure}[H]
	\centering
	\subfigure{\includegraphics[width=0.490\textwidth]{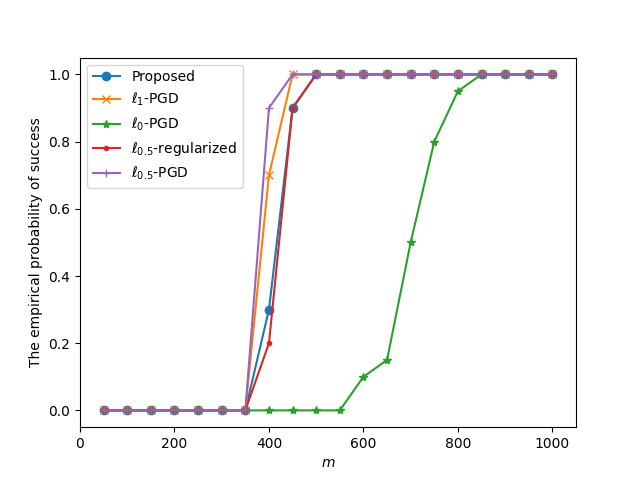}} 	\subfigure{\includegraphics[width=0.495\textwidth]{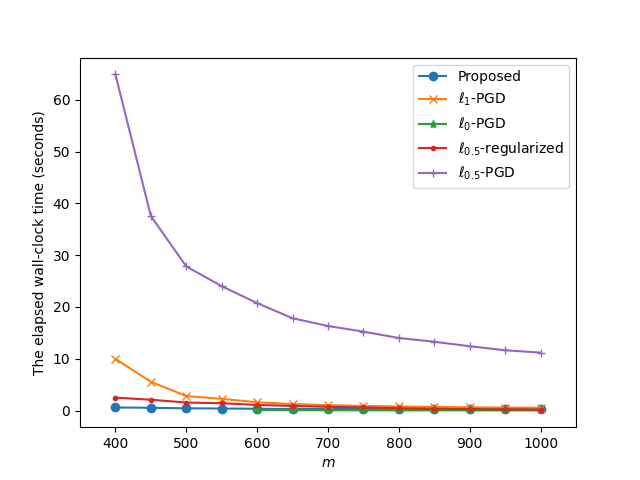}}\\
	\caption{Performance comparison of various considered algorithms. Left: The empirical probability of success versus $m$. Right: The elapsed wall-clock time versus $m$. }
	\label{fig.probability}
\end{figure}

We depict the empirical success probabilities by varying $m$ from $50$ to $1000$ in Figure \ref{fig.probability}, where the empirical probability of success is calculated as the ratio of successful recoveries to the total number of runs. The results for each $m$ are presented as an average over $20$ runs.
\vspace{-20pt}

According to Figure \ref{fig.probability}, we can first observe the effectiveness of all compared algorithms in addressing sparse recovery problems for $m \geq 550$ approximately. In terms of computational efficiency, the proposed algorithm stands out as the most time-efficient method. Other approaches exhibit similar computational times, with the exception of the $\ell_{0.5}$-PGD algorithm, which varies for different $m$.

\textbf{The nonconvex Cauchy loss.} In this experiment, we apply the proposed algorithm to the sparse recovery problem with a nonconvex loss function $f$. The results validate the generality of the proposed algorithm.

\begin{figure}[H]
	\centering
	\subfigure{\includegraphics[width=0.490\textwidth]{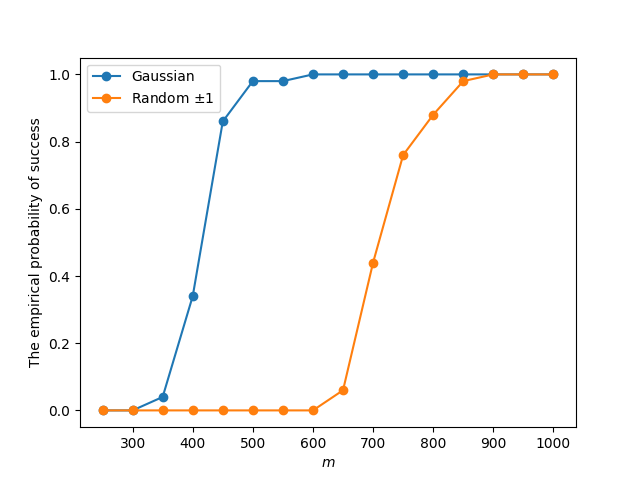}} 
	\subfigure{\includegraphics[width=0.495\textwidth]{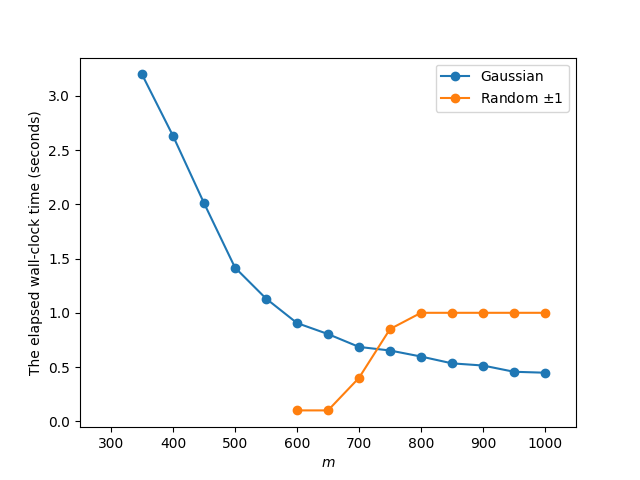}}\\
	\caption{Performance for the Cauchy loss. Left: The empirical probability of success versus $m$. Right: The elapsed wall-clock time versus $m$. }
	\label{fig.nonconvex}
\end{figure}

In addition to generating $\hat{\bm{x}}$ by randomly choosing values from $+1$ and $-1$, we conducted experiments with a modified $\hat{\bm{x}}$ distribution where nonzero entries follow a Gaussian distribution. This modification entails adjusting the initial step of the procedures as follows:
\begin{itemize}
	\item[(i)] Form $\hat{\bm{x}}$ by randomly selecting $n-s$ entries as zeros. Each nonzero entry is sampled from a standard normal distribution.
\end{itemize}
The initialization $\bm{x}^{0}$ and the stepsize $\beta = 1/\lambda_{\max}$ in the gradient projection step remain consistent with the convex squared loss test. The numerical results are illustrated in Figure \ref{fig.nonconvex}, and the presented results represent an average of over 20 runs. Specifically, we solely accounted for the average elapsed wall-clock time in instances where successful recovery occurred.
\vspace{-20pt}
\subsection{Image Reconstruction}
In this experiment, we employ the proposed algorithm to address \eqref{prob:sparse} in the context of real-world dataset image reconstruction. The underlying principle is rooted in the observation that natural images usually exhibit sparsity in wavelet domains, thereby reducing the required number of measurements for compressive imaging across diverse transformations.

Our experiments are conducted on the Set12 dataset \cite{zhang2017beyond}, which consists of several grayscale images sized at $256 \times 256$. We employ the discrete wavelet transform (DWT) for the sparse base, as described in \cite{liu2022improved}. Let $\hat{x}$ denote the original sparse bases, also known as wavelet coefficients in DWT, derived from the test images. Due to memory constraints, computing a measurement matrix of size $m \times 65536$ is impractical. Consequently, we treat the $256 \times 256$ wavelet coefficients as $256$ individual columns. 

\begin{figure}[htbp]
	\centering
	\includegraphics[width=4.5in]{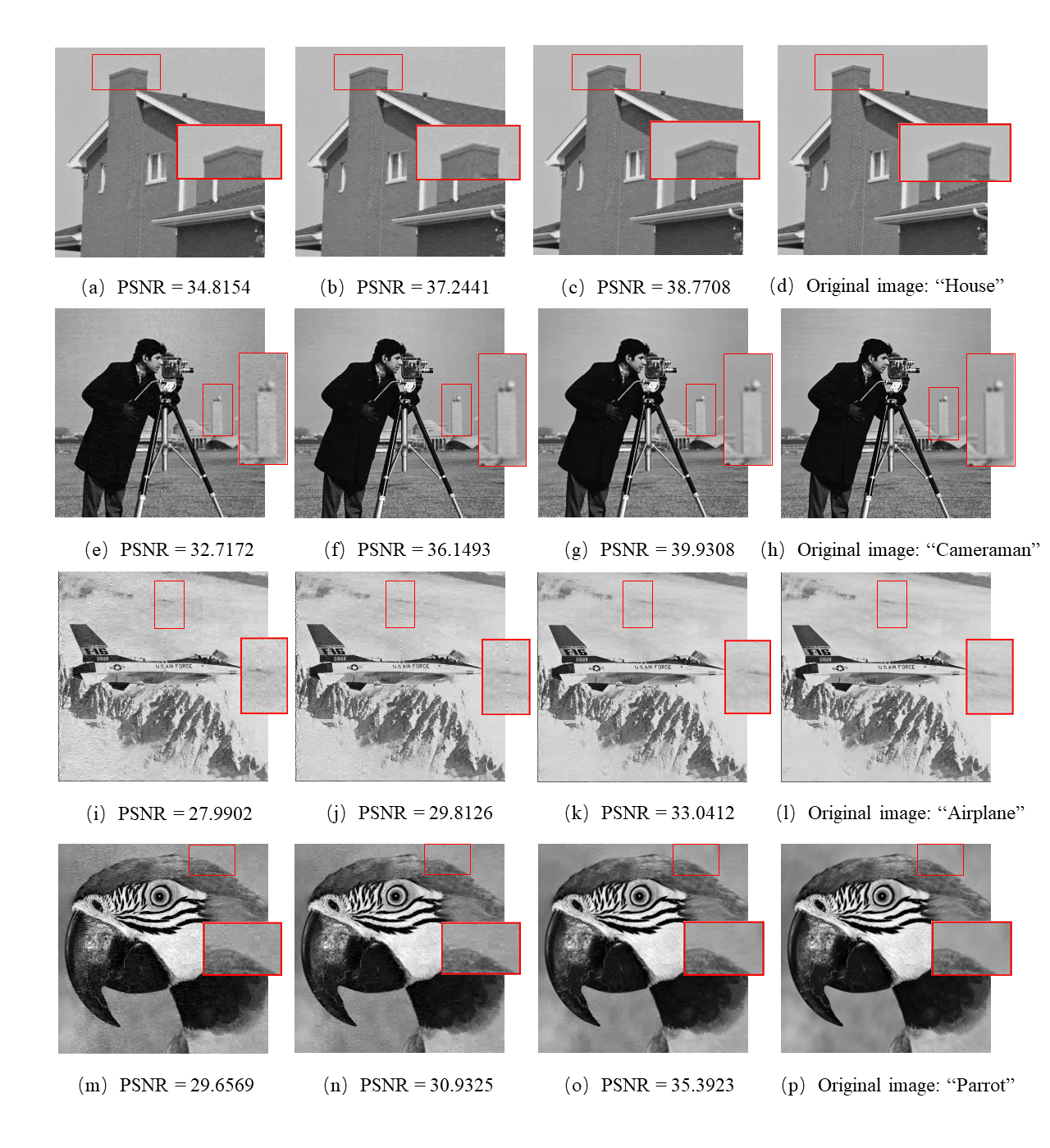}
	\caption{The first three columns are the reconstructed image with $p \in \{0.1,0.5,0.9\}$ and the last column corresponds to the original images.}
	\label{fig.image_recovery}
\end{figure}

The Peak Signal-to-Noise Ratio (PSNR) serves as the quantitative evaluation metric, defined as $\text{PSNR} = 20 \log_{10} \left( \frac{255}{\text{MSE}} \right)$, where $\text{MSE} = \frac{1}{mn} \sum_{i=1}^{m} \sum_{j=1}^{n} [I(i, j) - D(i, j)]^2$. Here, $I \in \mathbb{R}^{m \times n}$ and $D \in \mathbb{R}^{m \times n}$ represent the original and reconstructed images, respectively.
We mention that a higher PSNR value generally indicates greater similarity between the reconstructed image and the original image. We fix the number of measurements, $m$, at $200$, with varying values of $p$. For each image, we iterate through the following steps $10$ times:     
\begin{itemize}
	\item[(i)] Construct wavelet coefficients $\hat{\bm{x}}$ using DWT.
	\item[(ii)] Generate $ \bm{A} $ with entries drawn from a standard normal distribution.
	\item[(iii)]For each column in $\hat{\bm{x}}$, we first construct $\bm{b} = \bm{A}\hat{\bm{x}}$ and then address the corresponding optimization problems to derive the estimated wavelet coefficients.
	\item[(iv)] Combine the estimated wavelet coefficients of all columns to get $\bm{x}^{*}$.
	\item[(v)] Perform the inverse transformation on $\bm{x}^{*}$ to obtain the recovered vector.
\end{itemize}
The experimental setup of this experiment generally follows \S\ref{Sec_sparsesignal} involving the convex squared loss. The results are visually depicted in Figure \ref{fig.image_recovery}. 

\begin{figure}[htbp]
	\centering
	\begin{subfigure}{}
		\centering
		\includegraphics[width=4.5in]{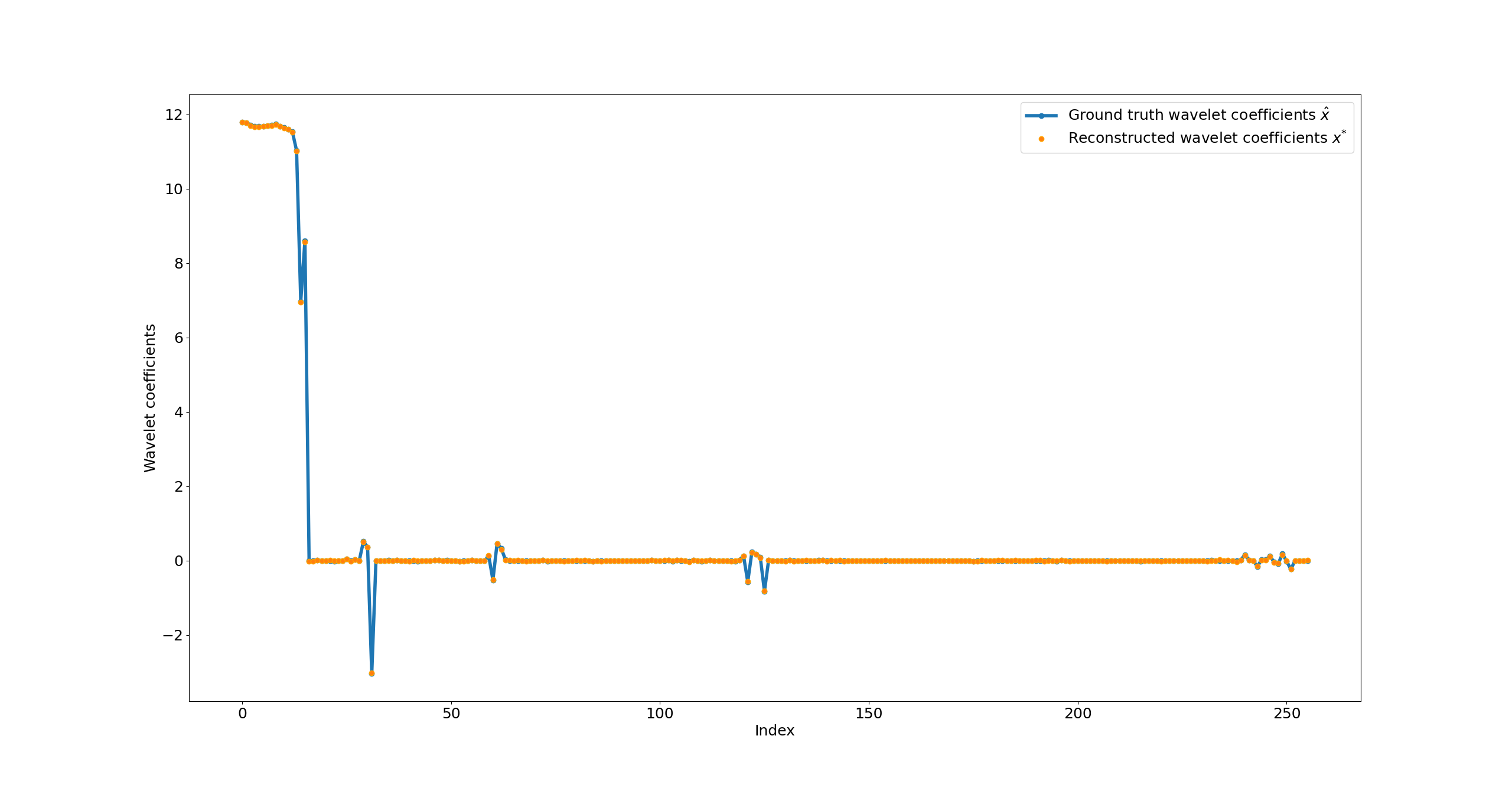}
		\caption{The 1st column of wavelet coefficients generated from image ``House".}
		\label{fig.wavelet1}
	\end{subfigure}%
	\begin{subfigure}{}
		\centering
		\includegraphics[width=4.5in]{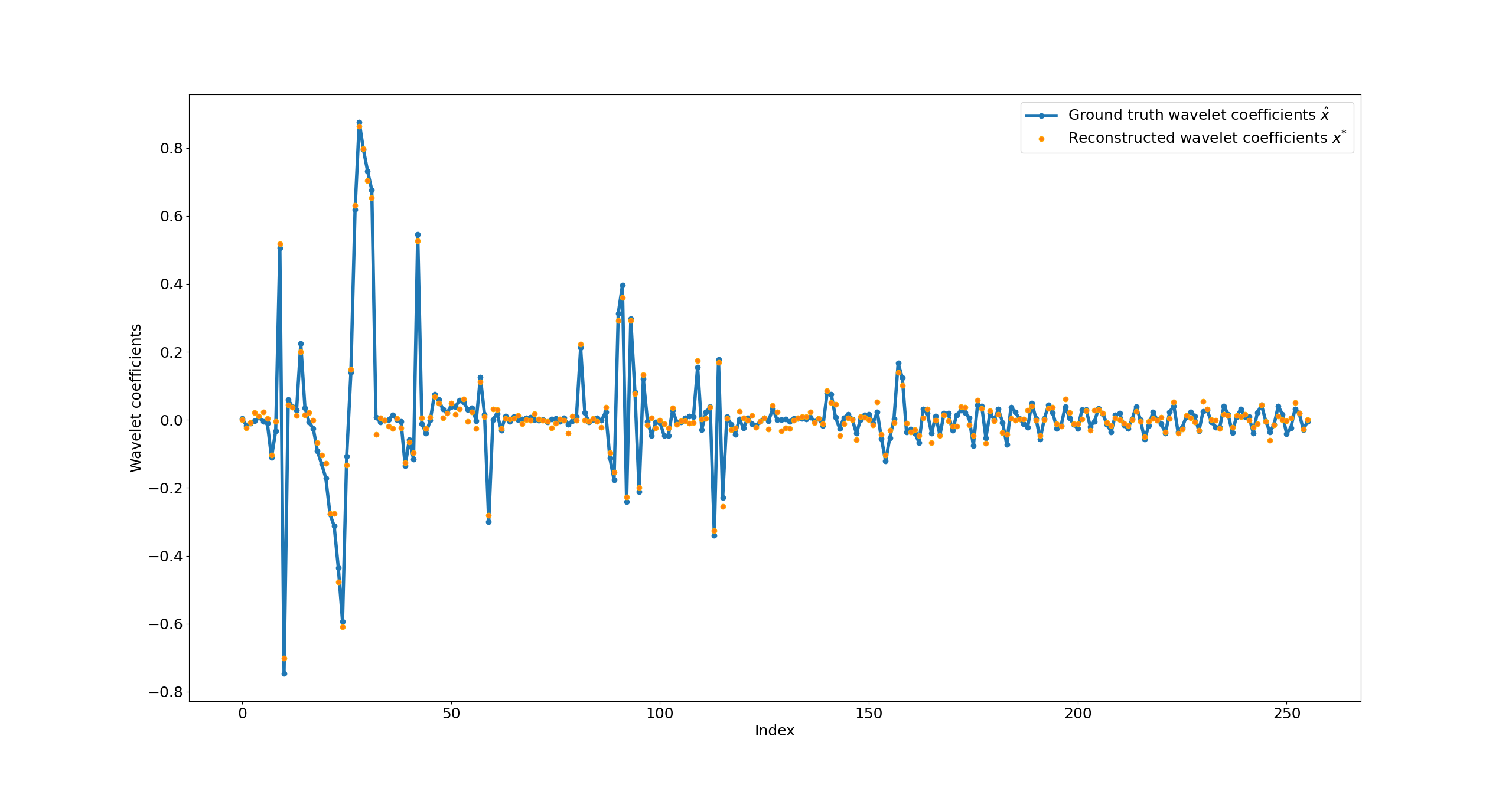}
		\caption{The 51st column of wavelet coefficients generated from image ``House".}
		\label{fig.wavelet2}
	\end{subfigure}%
\end{figure}

Figure \ref{fig.image_recovery} illustrates the effective recovery of original images through the proposed algorithm across various $p$ values. The presented results suggest that larger $p$ values appear to contribute to more efficient image recovery, primarily due to the sparsity level of the ground truth $\hat{\bm{x}}$. To clarify this correlation, we select two bases derived from the ``House" image, presenting the corresponding distributions of the reconstructed wavelet coefficients $\bm{x}^{*}$ and the ground truth wavelet coefficients of $\hat{\bm{x}}$ in Figure \ref{fig.wavelet1} and Figure \ref{fig.wavelet2}. In particular, Figure \ref{fig.wavelet2} reveals that the ground truth wavelet coefficients $\hat{\bm{x}}$ are not extremely sparse, as evidenced by numerous small yet nonzero wavelet coefficients. This phenomenon contributes to relatively large PSNR values for $p=0.1$, as observed in Figure \ref{fig.image_recovery}.

\section{Conclusion}\label{Con_sec}
In this paper, we have developed, analyzed, and implemented a novel hybrid first-order method to optimize a smooth function subject to nonconvex $\ell_{p}$ ball constraints. Notably, this approach can be readily extended to accommodate other level set constraints incorporating sparsity-promoting regularizers. The proposed method adaptively alternates between solving a Frank-Wolfe-type subproblem and a gradient projection subproblem, facilitating a simple yet efficient implementation. We have shown the global convergence properties of the proposed algorithm for nonconvex objectives and have established the local worst-case convergence rate of $O(1/\sqrt{k})$. Finally, we have evaluated the performance of the proposed algorithm through numerical studies focused on solving the large-scale $\ell_p$ ball projection problem, sparse signal recovery problem, and image reconstruction problem.

\begin{appendix}
	\section{Proofs of Propositions and Lemmas.}
	\subsection{Proof of Lemma \ref{Lem: Obj_decrease}} \label{app_C}
	\begin{proof}
		We first prove (i). Note that $\Vert \bm{x}\Vert_{p}\leq \gamma^{\frac{1}{p}}$ implies that $\Vert\bm{x}\Vert_1\leq\gamma^{\frac{1}{p}}$ and $\Vert\bm{x}\Vert_1 \geq \Vert\bm{x}\Vert_2$ for any $\bm{x} \in \mathbb{R}^{n}$. We have $\Vert\bm{x}\Vert_2 \leq \gamma^{\frac{1}{p}}$ provided that $\Vert \bm{x}\Vert_{p}\leq \gamma^{\frac{1}{p}}$.
		As a result, 
		\begin{equation}\label{eq: d_k inequality}
			\Vert\bm{d}^{k}\Vert_{2} = \Vert\bm{s}^{k} 	\!-\! \bm{x}^{k}\Vert_2\leq \Vert\bm{s}^{k} \Vert_2 \!+\! \Vert\bm{x}^{k} \Vert_2 \!\leq\! 2\gamma^{\frac{1}{p}},
		\end{equation}
		where the last inequality holds due to $\bm{s}^{k}$ and $\bm{x}^{k}$ are feasible for $\mathcal{B}_{\ell_p}$. 
		
		For $k \in \mathcal{T}_1$, since $ \alpha^k =\min
		(1, g_k/(L_{f}^k\Vert\bm{d}^k\Vert^2_2)) $ by Line \ref{line:localiniStep} of Algorithm \ref{alg.subprocedure_fw}, it then follows from the backtracking line-search stepsize that 
		\begin{equation}\label{eq: FW_reduction}
			\begin{aligned}
				f(\bm{x}^{k+1}) - f(\bm{x}^{k}) 
				&\leq -\bar{\alpha}^{k}g_{k} + \frac{(\bar{\alpha}^{k})^2 L_{k}^f}{2}\Vert \bm{d}^k \Vert_{2}^{2} \leq -\frac{(g_{k})^2}{2L_{f}^k\Vert \bm{d}^{k}\Vert_{2}^2} \\
				& \leq - \frac{g_{k}}{2} \min\left(\frac{g_k}{L_{f}^k\Vert \bm{d}^{k}\Vert_{2}^{2}},1 \right)\\
				& \leq - \frac{\Delta Q(\bm{s}^{k};\bm{x}^{k})}{2}\min\left(\frac{\Delta Q(\bm{s}^{k};\bm{x}^{k})}{4\max\left(\tau L_{f}, L_{-1}\right) \gamma^{2/p}},1 \right) \leq 0.
			\end{aligned}
		\end{equation}
		where the second inequality holds because $\bar{\alpha}^k$ minimizes the right-hand side of the unconstrained quadratic surrogate function in \eqref{eq: QuadrticApproximation} and the fourth inequality holds by making use of \cite[Proposition 2, Appendix C]{pedregosa2020linearly} and \eqref{eq: d_k inequality}.
		Similarly, for $k\in\mathcal{T}_2$, since $\alpha_{\textrm{bis}} \in (0,\bar{\alpha}^k)$, it holds that
		\begin{equation}\label{eq:Bisection_reduction}
			f(\bm{x}^{k+1}) - f(\bm{x}^{k}) \leq - \frac{\Delta Q^{2}(\bm{s}^{k};\bm{x}^{k})}{8\max\left(\tau L_{f}, L_{-1}\right) \gamma^{2/p}} \leq 0.
		\end{equation} 
		Rearranged, this leads to the desired inequality. 
		
		To prove the first inequality in (ii), we note that the proof follows the results of \cite[Lemma 3.6 and Theorem 3.7]{zhang2020inexact} and thus is omitted here. On the other hand, since the objective in \eqref{sub.proj} is strongly convex and $\bm{x}^{k+1}$ is the optimal solution of the $k$th subproblem, we have 
		\begin{equation}\label{eq: P_sub is decreasing}
			\begin{aligned}
				&\quad\ \Delta P(\bm{x}^{k+1};\bm{x}^{k})
				= P(\bm{x}^{k};\bm{x}^{k}) - P(\bm{x}^{k+1};\bm{x}^{k})\\ &=-\nabla f(\bm{x}^{k}_{\I^{k}})^{T}(\bm{x}^{k+1}_{\I^{k}} - \bm{x}^{k}_{\I^{k}}) - \frac{1}{2\beta}\Vert\bm{x}^{k+1}_{\I^{k}} - \bm{x}^{k}_{\I^{k}}\Vert_2^{2}\geq 0.
			\end{aligned}
		\end{equation}
		Recall $f$ has Lipschitz continuous gradient with Lipschit constant $L_f$ and by \eqref{eq: P_sub is decreasing}, we have
		\begin{equation}\label{eq: Proj_reduction}
			\begin{aligned}
				f(\bm{x}^{k+1})
				&\leq  f(\bm{x}^{k}) +  \nabla f(\bm{x}^{k})^{T}(\bm{x}^{k+1} - \bm{x}^{k}) + \frac{L_f}{2}\Vert\bm{x}^{k+1} - \bm{x}^{k}\Vert_2^{2}\\
				&\leq f(\bm{x}^{k}) + \nabla f(\bm{x}^{k})^{T}(\bm{x}^{k+1} - \bm{x}^{k}) \!+\! \frac{1}{2\beta}\Vert\bm{x}^{k+1} - \bm{x}^{k}\Vert_2^{2}\\
				&= f(\bm{x}^{k}) + \nabla f(\bm{x}^{k}_{\I^{k}})^{T}(\bm{x}^{k+1}_{\I^{k}} - \bm{x}^{k}_{\I^{k}}) + \frac{1}{2\beta}\Vert\bm{x}^{k+1}_{\I^{k}} - \bm{x}^{k}_{\I^{k}}\Vert_2^{2}\\
				&\leq f(\bm{x}^{k}) - \Delta P(\bm{x}^{k+1};\bm{x}^{k}),
			\end{aligned}
		\end{equation}
		where the equality holds due to the constraint imposed in \eqref{sub.proj}. Therefore, $\Delta P(\bm{x}^{k+1};\bm{x}^{k}) \leq \Delta f(\bm{x}^{k+1})$, as desired.
		
		We now prove (iii). Summing up both sides of \eqref{eq: FW_reduction}, \eqref{eq:Bisection_reduction} and \eqref{eq: Proj_reduction}  over $k$ gives
		\begin{equation}\label{eq: summing_up}
			\begin{aligned}
				&\quad\sum_{k\in \mathcal{T}_{1}} \min\left(\frac{\Delta Q^2(\bm{s}^{k};\bm{x}^{k})}{8\max\left(\tau L_{f}, L_{-1}\right) \gamma^{2/p}},\frac{\Delta Q(\bm{s}^{k};\bm{x}^{k})}{2} \right) 
				+ \sum_{k\in \mathcal{T}_{2}}\frac{\Delta Q^{2}(\bm{s}^{k};\bm{x}^{k})}{8\max\left(\tau L_{f}, L_{-1}\right) \gamma^{2/p}}\\
				&\quad\  + \sum_{k\in \mathcal{S}_{2}}  \Delta P(\bm{x}^{k+1};\bm{x}^{k})\\ 
				&\leq \sum_{k=0}^{t}\Delta f(\bm{x}^{k+1}) = f(\bm{x}^{0}) - f(\bm{x}^{t+1})
				\leq f(\bm{x}^{0}) - \underline{f} <+\infty.
			\end{aligned}
		\end{equation} 
		Letting $t \to +\infty$, it gives 
		\begin{equation*}
			\begin{aligned}
				\lim\limits_{\substack{k\in\mathcal{T}_1 \\ k\to +\infty}} \Delta Q(\bm{s}^{k};\bm{x}^{k}) = 0, \lim\limits_{\substack{k\in\mathcal{T}_2 \\ k\to +\infty}} \Delta Q(\bm{s}^{k};\bm{x}^{k}) = 0 \ \textrm{ and }
				\lim\limits_{\substack{k\in\mathcal{S}_2 \\ k\to +\infty}}\Delta P(\bm{x}^{k+1};\bm{x}^{k}) = 0.
			\end{aligned}
		\end{equation*}
		Consequently, $\lim\limits_{\substack{k\in\mathcal{S}_2 \\ k\to +\infty}}\Vert\bm{x}^{k+1} - \bm{x}^{k}\Vert_{2} = 0$.

		To prove (iv). If $\alpha_{\textrm{bis}}^{k}\to 0$, then we have $\lim\limits_{\substack{k\in\mathcal{T}_2 \\ k\to +\infty}} \Vert\bm{x}^{k+1} - \bm{x}^k\Vert_2^2 = 0$ since $\bm{x}^{k+1} = \bm{x}^k + \alpha^k_{\textrm{bis}}\bm{d}^k$ for $k\in \mathcal{T}_2$, and $\bm{d}^{k}$ is bounded by \eqref{eq: d_k inequality}. This completes the proof.
	 \end{proof}

	
	
	
	\subsection{Proof of Theorem \ref{Theo: FW_reformulation}}\label{app_Sub} 
	The techniques of the proof of Theorem \ref{Theo: FW_reformulation} are mainly from the geometric approach perspective. The following lemma first helps characterize some basic properties of the global minimizers of \eqref{eq: FW_linearization}.
	\begin{lemma}\label{Lem: basic_properties of FW_sub}
		Suppose $\nabla f(\bm{x}^{k}) \neq \bm{0}$. Let $\bm{s}^{k}$ be a global optimal solution of subproblem \eqref{eq: FW_linearization}. Then the following statements hold:
		\begin{itemize}
			\item[(i)] $\text{sgn}(\bm{s}^{k}) = -\text{sgn}(\nabla f(\bm{x}^{k}))$, i.e., the corresponding entries of $\nabla f(\bm{x}^{k})$ and $\bm{s}^{k}$ have the opposite signs.
			
			\item[(ii)] $\Vert\bm{s}^{k}\Vert_{p}^{p} = \gamma$, i.e., $\bm{s}^{k}$ lies on the boundary of the $\ell_{p}$ ball. Indeed, $\bm{s}^{k}$ is a vertex of $\mathcal{B}_{\ell_p}$. Consequently, $\bm{s}^{k} = \argmin\limits_{\bm{s} \in \textrm{conv}(\mathcal{B}_{\ell_p})} Q(\bm{s};\bm{x}^k)$.  
		\end{itemize}
	\end{lemma}
	\begin{proof}		
		(i) Trivial. For (ii), we prove this by contradiction. Suppose this is not true, i.e., $\Vert\bm{s}^{k}\Vert_{p}^{p} = \sum\limits_{i\in\I(\bm{s}^{k})}\vert s_i^k\vert^{p} < \gamma$. Note that $\nabla_i f(\bm{x}^{k}) \neq 0$ for any $i \in \mathcal{I}(\bm{s}^{k})$, we choose an arbitrary $j \in \mathcal{I}(\bm{s}^{k})$.
		Without loss of generality, let $\tilde{\bm{s}}$ be a solution candidate such that $\tilde{s}_i = s^k_i$ for any $ i \in \mathcal{I}(\bm{s}^{k})\backslash \{j\} $ satisfying $ \vert\tilde{s}_j\vert > \vert s^{k}_j\vert$, $\Vert \tilde{\bm{s}} \Vert_{p}^{p} \leq \gamma$ and $\text{sgn}(\tilde{s}_j)$ = -$\text{sgn}(\nabla_j f(\bm{x}^{k}))$. We thus have $Q(\tilde{\bm{s}};\bm{x}^k) = \la \nabla f(\bm{x}^{k}), \tilde{\bm{s}}\ra = \sum\limits_{i\in\I(\bm{s}^{k})\backslash \{j\}}\nabla_i f(\bm{x}^{k})s_i^k + \nabla_j f(\bm{x}^{k})\tilde{s}_j = \sum\limits_{i\in\I(\bm{s}^{k})\backslash \{j\}}\nabla_i f(\bm{x}^{k})s_i^k - \vert \nabla_j f(\bm{x}^{k}) \vert \vert \tilde{s}_j\vert < \sum\limits_{i\in\I(\bm{s}^{k})\backslash \{j\}}\nabla_i f(\bm{x}^{k})s_i^k - \vert\nabla_j f(\bm{x}^{k})\vert \vert s_j^k\vert = \sum\limits_{i\in\I(\bm{s}^{k})\backslash \{j\}}\nabla_i f(\bm{x}^{k})s_i^k + \nabla_j f(\bm{x}^{k}) s_j^k = Q(\bm{s}^{k};\bm{x}^{k})$. This contradicts the optimality of $\bm{s}^{k}$. 
		
		Moreover, it follows from the linearity of the objective function $Q(\bm{s};\bm{x}^{k})$ and $\Vert\bm{s}^{k}\Vert_{p}^{p} = \gamma$ that $\bm{s}^{k}$ resides at the vertex of $\mathcal{B}_{\ell_p}$. Hence, this, together with \cite[Theorem 2.7 and Theorem 2.9]{bertsimas1997introduction} that $\bm{s}^{k}$ can be obtained by minimizing $Q(\bm{s};\bm{x}^k)$ over the convex hull of the nonconvex $\mathcal{B}_{\ell_p}$. This completes the proof.
	 \end{proof}
	
	We now prove Theorem \ref{Theo: FW_reformulation} in the following.
	\begin{proof}
		By Lemma \ref{Lem: basic_properties of FW_sub} and the minimization nature of \eqref{eq: FW_linearization}, we know that $\text{sgn}(s_{i_{\max}}^k) = -\text{sgn}(\nabla_{i} f(\bm{x}^{k}))$ and the rest entries of $\bm{s}^k$ are $0$. It follows that $\sum_{i=1}^{n}\vert s_i^k\vert^p = \vert s_{i_{\max}}^k \vert^{p} = \gamma$. Thanks to Assumption \ref{BasicAssum}(iii), we have that $\phi^{-1}(\vert s_{i_{\max}}^k \vert^{p}) = \vert s_{i_{\max}}^k \vert = \phi^{-1}(\gamma)$ (applied with $\phi(\cdot) = \vert \cdot \vert^{p}$). Therefore, we achieve the desired conclusion.
	 \end{proof}
	
	
	\section{An Alternating Projection Algorithm onto the Weighted $\ell_{1}$ Ball}\label{Appendix_C}
	
	We present an algorithm for computing the projection onto the weighted $\ell_{1}$ ball, which is an extension of the projection algorithm onto unweighted $\ell_{1}$ ball \cite{condat2016fast,zhang2020inexact}. By the symmetry of the weighted $\ell_{1}$ ball, we only have to consider the following projection problem with $\bm{y}\in \mathbb{R}^{n}_{+}$.
	
	\begin{equation}\label{eq: sub5.1}
		\begin{aligned}
			\min_{\bm{x}} \ \ \frac{1}{2}\Vert\bm{x}-\bm{y}\Vert_2^2, \quad \textrm{s.t.} \ \ \bm{x} \in \Delta \subset \mathbb{R}^{n}.
		\end{aligned}
	\end{equation}
	Here $ \Delta:=\{\bm{x}\in\mathbb{R}^{n}\mid\sum_{i=1}^nw_ix_{i} = \varrho \textrm{ and } x_{i} \geq 0, \,\, \forall i\in \left[n\right]\}$ denotes the weighted simplex. The projection onto the weighted simplex can be expressed as 
	$$\mathcal{P}_{\Delta}(\bm{y}) = \max\{\bm{0}, \bm{y} - \lambda\bm{w}\},$$
	where $\lambda > 0$ satisfies $\sum_{i \in [n]}w_i\max\{0,y_i-\lambda w_i\} = \varrho$. As analyzed in \cite{zhang2020inexact}, the feasible set of \eqref{eq: sub5.1} can be expressed as the intersection of a hyperplane and  $\mathbb{R}^{n}_{+}$, so that \eqref{eq: sub5.1} can be recast as
	
	\begin{equation}\label{eq: alternating}
		\begin{aligned}
			\min_{\bm{x}} \ \ \frac{1}{2}\Vert\bm{x}-\bm{y}\Vert_2^2, \quad
			\textrm{s.t.} \ \ \bm{x}\in \mathcal{H}\cap\mathbb{R}^{n}_{+},
		\end{aligned}
	\end{equation}
	where $\mathcal{H} = \{\bm{x}\in \mathbb{R}^{n} \mid \sum_{i=1}^{n} w_ix_i = \varrho\}$. The projection of $\bm{y}$ onto $\mathcal{H}$ is given by
	\begin{equation}\label{eq: projection_intersection1}
		\mathcal{P}_{\mathcal{H}}(\bm{y}) = \bm{y}- \frac{\langle\bm{w},\bm{y}\rangle-\varrho}{\Vert\bm{w}\Vert_2^2}\bm{w}.
	\end{equation}
	One can easily extend the results in \cite[Lemma 4.2]{zhang2020inexact} to the following results.
	\begin{lemma}\label{lem: active-set_property}
		Let $\bm{x} = \mathcal{P}_{\mathcal{H}}(\bm{y})$. It holds that
		\begin{itemize}
			\item[($i$)] if $\bm{x} \geq \bm{0}$, then $\bm{x} = \mathcal{P}_{\mathcal{H}}(\bm{y}) = \mathcal{P}_{\Delta}(\bm{y})$.

			\item[$(ii)$] if $x_i \leq 0$, then $[\mathcal{P}_{\Delta}(\bm{y})]_i = 0, \forall i\in [n]$.
		\end{itemize}
	\end{lemma}
	
	From Lemma \ref{lem: active-set_property}, after obtaining $\mathcal{P}_{\Delta}(\bm{y})$, we can discard its nonpositive components. In this way, we can form another hyperplane in the subspace consisting of the nonzeros in  $\mathcal{P}_{\Delta}(\bm{y})$ and repeat the same procedure until there is no negative component in the projection onto the hyperplane. This is basically the same algorithm for the $\ell_{1}$-ball projection presented in \cite{zhang2020inexact}. It is easy to see the stated above procedure terminates in finite steps. In practice, $O(n)$ complexity can be often witnessed, e.g., see \cite{condat2016fast} and the references therein.
	
	The projection algorithm is summarized in Algorithm \ref{alg:l2}.
	
	\begin{algorithm}[H]
		\caption{An Alternating Projection Algorithm onto the Weighted $\ell_1$-Ball}
		\label{alg:l2}
		\begin{algorithmic}[1]
			\State \textbf{Input:} $\varrho > 0$ and $\bm{y}\in\mathbb{R}^{n}_{+}$.
			\State \textbf{Initialization:} $\bm{z}^{0} = \bm{y}$, set $j = 0$, $\mathcal{M}^{0} = \N^{0} = [n]$.
			\While{$\N^{j}\neq \emptyset$}
			\State Compute $\bm{x}^{j+1}_{\M^{j}} = 
			\bm{z}^{j+1} - \tfrac{\langle\bm{w}_{\M^{j}},\bm{y}_{\M^{j}}\rangle-\varrho}{\Vert\bm{w}_{\M^{j}}\Vert_2^2} \bm{w}_{\M^{j}}$ and $\bm{x}^{j+1}_{[n]\backslash \M^{j}} = \bm{0}$.
			\State Update $\M^{j+1} = \{i\mid x_{i}^{j+1} > 0, i\in[n]\}$ and  $\mathcal{N}^{j+1} = \{i\mid x_{i}^{j+1} < 0, i\in[n]\}$.
			\State Set $\bm{z}^{j+1} = \bm{x}^{j+1}_{\M^{j+1}}$, $\bm{w} = \bm{w}_{\M^{j+1}}$.
			\State Set $j\gets j+1$.
			\EndWhile
			\State \textbf{Output:}~$\bm{x}^{k+1}$.
		\end{algorithmic}
	\end{algorithm}
\end{appendix}

\bibliographystyle{plain}
\bibliography{references}

\end{document}